\documentclass{article}

\usepackage{arxiv}

\usepackage[utf8]{inputenc} 
\usepackage[T1]{fontenc}    
\usepackage{url}            
\usepackage{booktabs}       
\usepackage{amsfonts}       
\usepackage{nicefrac}       
\usepackage{microtype}      
\usepackage{lipsum}		
\usepackage{subfig}
\usepackage{graphicx}
\usepackage{doi}

\usepackage{amsmath,amssymb,amsfonts}
\usepackage{graphicx}
\usepackage{textcomp}
\def\BibTeX{{\rm B\kern-.05em{\sc i\kern-.025em b}\kern-.08em
    T\kern-.1667em\lower.7ex\hbox{E}\kern-.125emX}}

\usepackage{bbm}
\usepackage{multirow}
\usepackage{mathtools}
\usepackage{cuted}
\usepackage{algpseudocode,algorithm}

\DeclareMathAlphabet\mathbfcal{OMS}{cmsy}{b}{n}

\newcommand\Hyp[3]{\left(\begin{array}{@{\,}c@{\,}}#1\\#2\end{array}\middle\lvert\;#3\right)}
\newcommand\drp[2]{\mathbfcal{R}^{\left(\substack{\alpha,\beta\\a,b}\right)}_{#1}\!\!\left(#2;N\right)}

\newcommand\drpn[2]{\hat{\mathbfcal{R}}^{\left(\substack{\alpha,\beta\\a,b}\right)}_{#1}\!\!\left(#2;N\right)}
\newcommand\drpnn[3]{\hat{\mathbfcal{R}}^{\left(\substack{\alpha_{#1},\beta_{#1}\\a_{#1},b_{#1}}\right)}_{#2}\!\!\left(#3\right)}

\newcommand\drpo[2]{\hat{\mathbfcal{R}}_{#1}\!\left(#2\right)}
\newcommand\drps[2]{\hat{\mathbfcal{R}}^{\left(\substack{0, 0\\0,b}\right)}_{#1}\!\!\left(#2\right)}
\newcommand\abs[1]{\left\lvert#1\right\rvert}

\title{Accelerated and Improved Stabilization for High Order Moments of Racah Polynomials}

\author{ \href{https://orcid.org/0000-0002-4121-0843}{\includegraphics[scale=0.06]{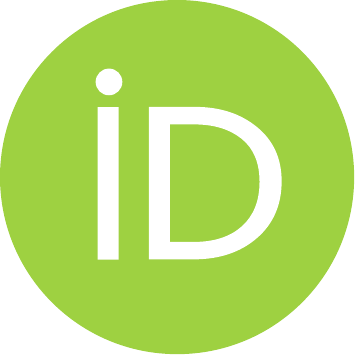}\hspace{1mm}Basheera M. Mahmmod}
\\	
    Department of Computer Engineering\\
	University of Baghdad\\
	Baghdad, 10071, Iraq \\
	\texttt{basheera.m@coeng.uobaghdad.edu.iq} \\
    \And
	\href{https://orcid.org/0000-0002-6439-0082}{\includegraphics[scale=0.06]{orcid.pdf}\hspace{1mm}Sadiq H.  {Abdulhussain}} \\
    Department of Computer Engineering\\
	University of Baghdad\\
	Baghdad, 10071, Iraq \\
    \texttt{sadiqhabeeb@coeng.uobaghdad.edu.iq} \\
	\And
	\href{https://orcid.org/0000-0001-6737-0565}{\includegraphics[scale=0.06]{orcid.pdf}\hspace{1mm}Tom\'{a}\v{s} {Suk}} \\
	Czech Academy of Sciences \\
	Institute of Information Theory and Automation\\
	Pod vod\'{a}renskou v\v{e}\v{z}\'{\i} 4, Praha 8, 182\,08,  Czech Republic \\
	\texttt{suk@utia.cas.cz} \\
}

\begin{document}
\maketitle

\begin{abstract}
One of the most effective orthogonal moments, discrete Racah polynomials (DRPs) and their moments are used in many disciplines of sciences, including image processing, and computer vision.
Moments are the projections of a signal on the polynomial basis functions. Racah polynomials were introduced by Wilson and modified by Zhu for image processing and they are orthogonal on a discrete set of samples. However, when the moment order is high, they experience the issue of numerical instability. In this paper, we propose a new algorithm for the computation of DRPs coefficients called Improved Stabilization (ImSt). In the proposed algorithm, {the DRP plane is partitioned into four parts, which are asymmetric because they rely on the values of the polynomial size and the DRP parameters.} The logarithmic gamma function is utilized to compute the initial values, which empower the computation of the initial value for a wide range of DRP parameter values as well as large size of the polynomials. In addition, a new formula is used to compute the values of the initial sets based on the initial value. Moreover, we optimized the use of the stabilizing condition in specific parts of the algorithm. ImSt works for wider range of parameters until higher degree than the current algorithms. We compare it with the other methods in a number of experiments.
\end{abstract}

\keywords{Racah polynomials \and Orthogonal moments \and Recurrence algorithm \and Stabilizing condition}

\section{Introduction}
Moment can be understood as the projection of a signal to a polynomial basis. The moments are widely used as features for recognition of images and various image-like data. The moments can be divided to non-orthogonal and orthogonal. The non-orthogonal geometric and complex moments have advantage in easier construction of invariants to various geometric and radiometric transformations, e.g. rotation \cite{hosny2018new}, \cite{leonid1}, affine transformation \cite{Hickman}, convolution with symmetric filter \cite{FluBolZit}, \cite{boldys08} etc. 
On the other hand, they are very correlated each other, what leads to precision loss in lower orders than that of the orthogonal moments (the order equals the degree of the polynomial).

That is why we use orthogonal polynomials. They can be further divided into continuous and discrete. The relation of orthogonality of the continuous polynomials is based on integral over some interval, an example can be Fourier Mellin moments \cite{Singh_orthogonal_Fourier_Mellin}. When we compute the continuous moment from a digital image that is only defined in discrete pixels, we obtain the value with some error caused by the approximate computation of the definition integrals. Therefore the polynomials with discrete orthogonality that is based on the sum over some finite set of discrete samples are intensively studied.

Different types of discrete orthogonal polynomials have been derived over the ages. Here, we mention only that with significance for image processing. Besides his famous continuous polynomials, Chebyshev published also discrete ones. Mukundan derived efficient algorithm for computation of the discrete Chebyshev polynomials \cite{mukundan3}. Krawtchouk polynomials have parameter $p\in\langle 0, 1 \rangle$. It moves the zeros over the image and we can use it for adjustment of the region of interest. Efficient algorithm for their computation can be found in \cite{Sadiq_Krawtchouk} or \cite{abdulhussain2019new}, non-traditional way of computation by filters was published in \cite{Honarvar_Flusser_Krawtchouk}. The generalization of the Krawtchouk polynomials are Meixner polynomials; the efficient algorithm is in \cite{Sadiq_Meixner}.

Other group of discrete orthogonal polynomials contains e.g. Hahn polynomials. They can be computed by the algorithm from \cite{Sadiq_Hahn}. Dual Hahn polynomials were derived by swapping coordinate and order of the Hahn polynomials. The result is the non-uniform lattice $x(s)=s(s+\gamma+\delta+1)$, see \cite{Koekoek}. It is difficult to use in image processing, therefore Zhu et al. \cite{Zhu:dual_hahn} slightly changed the definition and used the index $s$ as coordinate in the digital image. Efficient algorithm can be found in \cite{Daoui_dual_Hahn}.

The Racah polynomials were first published by Wilson in \cite{Wilson_Racah} and named after physicist and mathematician Giulio Racah. They have the similar non-uniform lattice $x(s)=s(s+\gamma+\delta+1)$, as the dual Hahn, see \cite{Koekoek}. Zhu et al. \cite{Zhu_Racah} solved it also similarly. The Racah moments were used in skeletonization of craft images \cite{Fardousse_Racah}, Chinese character recognition \cite{Wu_Chinese_characters_Racah}, handwritten digit recognition \cite{Salouan_arabic}, and face recognition \cite{Ananth_Bharathi_face}.

In this paper, we propose an efficient algorithm for computation of the Racah polynomials.
The paper is organized as follows. Sec.~\ref{Prelim} is summary of definitions and state-of-the-art algorithms, our proposed method is in Sec.~\ref{Method}, we show its properties in numerical experiments in Sec.~\ref{Exper} and Sec.~\ref{Conc} concludes the paper.

recognition \cite{Salouan_arabic}, and face recognition \cite{Ananth_Bharathi_face}.

In this paper, we propose an efficient algorithm for computation of the Racah polynomials.
The paper is organized as follows. Sec.~\ref{Prelim} is summary of definitions and state-of-the-art algorithms, our proposed method is in Sec.~\ref{Method}, we show its properties in numerical experiments in Sec.~\ref{Exper} and Sec.~\ref{Conc} concludes the paper.

\section{Preliminaries and Related Work}
\label{Prelim}
In this section, the mathematical definitions and fundamentals of the discrete Racach polynomial (DRP) and their moments are presented. The current methods of their computation are summarized.

\subsection{The mathematical definition of DRPs}
The original Wilson's definition \cite{Wilson_Racah} is
\begin{equation}
\label{Eq_Wilson}
\mathbfcal{R}^{\left(\substack{\alpha,\beta\\ \gamma,\delta}\right)}_{n}\!\!\left(\lambda(x);N\right)
\!=\!{}_4F_3\Hyp{-n,n+\alpha+\beta+1,-x,x+\gamma+\delta+1}{\alpha+1,\beta+\delta+1,\gamma+1}{1},
\end{equation}
where $ {}_4F_3(\cdot) $ is the hypergeometric series. It is defined 
\begin{equation}
\label{Eq_hyp}
{}_4F_3\Hyp{a,b,c,d}{e,f,g}{z} = \sum\limits_{k=0}^{\infty} \frac{(a)_k \, (b)_k \, (c)_k\, (d)_k}{(e)_k \, (f)_k \, \, (g)_k}\cdot \frac{\left(z\right)^k}{k!}\ .
\end{equation}
The symbol $ (\cdot)_m $ is the Pochhammer symbol defined as
\begin{equation}
\label{Eq_pochhamer}
	(a)_m = a(a+1)(a+2)\cdots(a+m-1)\ .
\end{equation}
Zhu et al. in \cite{Zhu_Racah} introduced a new variable $s$ and defined $x=s(s+1)$.
Then the $ n $th order of the DRPs $ \drp{n}{s} $ are given by
\begin{align}
\label{Eq1}
	\drp{n}{s}&=\frac{1}{n!}(a+b+\alpha+1)_n(\beta+1)_n(a-b+1)_n \times {}_4F_3\Hyp{-n,a-s,a+s+1,\alpha+\beta+n+1}{\beta+1,a+b+\alpha+1,a-b+1}{1},
\end{align}
where $a>-1/2$, $b>a$, $b-a=N$ must be integer, $\alpha>-1$, $\beta>-1$, and $\beta<2a+1$.
    
The DRPs satisfy the condition of orthogonality
\begin{equation}
\label{Eq_oc}
\sum\limits_{s=a}^{b-1}\drp{n}{s} \drp{m}{s} \rho(s) \Delta x\left(s-\frac{1}{2}\right)= d_n^2 \delta_{nm}\ ,
\end{equation}
where $ \delta_{nm} $ is the Kronecker delta, $\Delta x\left(s-\frac{1}{2}\right)$ is the difference of the $x$ shifted by a half, i.e. $\Delta x\left(s-\frac{1}{2}\right)=\left(s+\frac{1}{2}\right)\left(s+\frac{3}{2}\right)-\left(s-\frac{1}{2}\right)\left(s+\frac{1}{2}\right)=(2s+1) $, $ \rho $ is the weight function of DRP
\begin{equation}
\label{weight}
	\rho(s)\!=\!\frac{\Gamma(a\!+\!s\!+\!1)\Gamma(b\!+\!s\!+\!\alpha\!+\!1)\Gamma(b\!+\!\alpha\!-\!s)\Gamma(s\!-\!a\!+\!\beta\!+\!1)}{\Gamma(b\!+\!s\!+\!1)\Gamma(b\!-\!s)\Gamma(s\!-\!a\!+\!1)\Gamma(a\!-\!\beta\!+\!s\!+\!1)}
\end{equation}
and $ d_n^2 $ is the norm function of DRP
\begin{equation}
\label{norm}
d_n^2\!=\!\frac{\Gamma(\alpha\!+\!n\!+\!1)\Gamma(\beta\!+\!n\!+\!1)\Gamma(a\!+\!b\!+\!\alpha\!+\!n\!+\!1)\Gamma(b\!-\!a\!+\!\alpha\!+\!\beta\!+\!n\!+\!1)}{(\alpha\!+\!\beta\!+\!2n\!+\!1)\Gamma(n\!+\!1)\Gamma(b\!-\!a\!-\!n)\Gamma(\alpha\!+\!\beta\!+\!n\!+\!1)\Gamma(a\!+\!b\!-\!n\!-\!\beta)}
\end{equation}
The $ n $th degree of the weighted DRP is given by
\begin{equation}\label{Eq_13}
	\drpn{n}{s} = \drp{n}{s}\sqrt{\frac{\rho_(s)}{d_n^2}\cdot\Delta x\left(s-\frac{1}{2}\right)} .
\end{equation}

\subsection{The state-of-the-art algorithms}
We can find significant algorithms of two authors in the literature, the original Zhu's paper and Daoui's approach.
For convenience, we will use the simplified notation $\drpn{n}{s}=\drpo{n}{s}$ with $b=a+N$.

\subsubsection{Zhu's algorithms}
Zhu et al. in \cite{Zhu_Racah} published two algorithms for Racah polynomial computation: recurrence over the order $n$ and recurrence over the index $s$.
The recurrence formula of weighted Racah polynomials over the order $n$ is
\begin{equation}
\label{Racah_recurrence_n}
\drpo{n+1}{s}=\left(B\frac{d_{n}}{d_{n+1}}\drpo{n}{s}-C\frac{d_{n-1}}{d_{n+1}}\drpo{n-1}{s}\right)/A
\end{equation}
with initial conditions
\begin{equation}
\label{Racah_recurrence_n_init}
\begin{array}{l}
\drpo{0}{s}\!=\!\sqrt{\frac{\varrho(s)}{d_n^2}\scriptstyle{(2s+1)}},\\
\drpo{1}{s}\!=\!-\sqrt{\frac{\varrho(s)}{d_n^2}\scriptstyle{(2s\!+\!1)}}\times\left(\frac{\varrho(s\!+\!1)(s\!+\!1\!-\!a)(s\!+\!1\!+\!b)(s\!+\!1\!+\!a\!-\!\beta)(b\!+\!\alpha\!-\!s\!-\!1)}{\varrho(s)(2s\!+\!1)}\!-\!\frac{\varrho(s)(s\!-\!a)(s\!+\!b)(s\!+\!a\!-\!\beta)(b\!+\!\alpha\!-\!s)}{\varrho(s)(2s\!+\!1)}\right),
\end{array}
\end{equation}
where
\begin{equation}
\label{Racah_recurrence_n_coeff}
\begin{array}{l}
A\!=\frac{(n+1)(\alpha+\beta+n)}{(\alpha+\beta+2n+1)(\alpha+\beta+2n+2)},\\[0.6ex]
B\!=\scriptstyle{s(s+1)-}\frac{\scriptstyle{a^2\!+b^2\!+(a-\beta)^2\!+(b+\alpha)^2}}{\scriptstyle{4}}+\frac{\scriptstyle{(\alpha+\beta+2n)(\alpha+\beta+2n+2)}}{\scriptstyle{8}}-
-\frac{\scriptstyle{(\beta^2\!-\alpha^2)[(2b+\alpha)^2\!-(2a-\beta)^2]}}{\scriptstyle{8(\alpha+\beta+2n)(\alpha+\beta+2n+2)}},\\[0.6ex]
C\!=\frac{(\alpha+n)(\beta+n)}{(\alpha+\beta+2n)(\alpha+\beta+2n+1)} \!\left[\left(\scriptstyle{a+b+}\frac{\scriptstyle{\alpha-\beta}}{\scriptstyle{2}}\right)^2\!\scriptstyle{-}\left(\scriptstyle{n}+\frac{\scriptstyle{\alpha+\beta}}{\scriptstyle{2}}\right)^2\right]\!
\!\left[\left(\scriptstyle{b-a+}\frac{\scriptstyle{\alpha+\beta}}{\scriptstyle{2}}\right)^2\!\scriptstyle{-}\left(\scriptstyle{n}+\frac{\scriptstyle{\alpha+\beta}}{\scriptstyle{2}}\right)^2\right]\!.
\end{array}
\end{equation}

The second algorithm is recurrence over the index $s$
\begin{equation}
\label{Racah_recurrence_s}
\begin{array}{l}
\drpo{n}{s} = \frac{(2s-1)[\sigma(s-1)+(s-1)\tau(s-1)-2\lambda s(s-1)]}{(s-1)[\sigma(s-1)+(2s-1)\tau(s-1)]}\sqrt{\frac{\rho(s)(2s+1)}{\rho(s-1)(2s-1)}}\drpo{n}{s-1}-\\[1ex]
\hspace{4.2em}\frac{2\sigma(s-1)}{(s-1)[\sigma(s-1)+(2s-1)\tau(s-1)]}\sqrt{\frac{\rho(s)(2s+1)}{\rho(s-2)(2s-3)}}\drpo{n}{s-2},
\end{array}
\end{equation}
where
\begin{equation}
\label{Racah_recurrence_s_func}
\begin{array}{l}
\sigma(s)=(s+a-\beta(b+\alpha-s))(s-a)(s+b)\\
\tau(s)=a(\alpha{+}1)(a{-}\beta){+}b(b{+}\alpha)(\beta{+}1){-}(\alpha{+}1)(\beta{+}1){-}s(s{+}1)(\alpha{+}\beta{+}2)\\
\lambda=n(n+1+\alpha+\beta).
\end{array}
\end{equation}
The declared initial values are
\begin{equation}
\label{Racah_recurrence_n_init_cond}
\begin{array}{l}
\drpo{n}{a}=\frac{(-1)^n}{(n!)^2}(a+1)_n(\beta-a+1)_n(b+\alpha+1)_n(b-n)_n \sqrt{\frac{\rho(0)}{d_n^2}},\\
\drpo{n}{a+1}=\frac{2}{(n+2)(n+1)}\left[\frac{\rho_n(1)}{\rho_n(0)}-\frac{n(n+1)}{2}\right] \sqrt{\frac{3\rho(1)}{\rho(0)}} \drpo{0}{s},
\end{array}
\end{equation}
where
\begin{equation}
\label{weightn}
\begin{array}{l}
\rho_n(s)=\frac{\Gamma(a+s+n+1)\Gamma(s-a+\beta+n+1)\Gamma(b+\alpha-s)\Gamma(b+\alpha+s+n+1)}{\Gamma(a-\beta+s+1)\Gamma(s-a+1)\Gamma(b-s-n)\Gamma(b+s+1)}.
\end{array}
\end{equation}
These initial conditions does not work; to resolve this issue, either the recurrence over $n$ for $s=a$ and $s=a+1$ is used or one of the following algorithms can be used.

\subsubsection{Daoui's algorithms}
Daoui et al. in \cite{daoui2022} proposed stabler algorithm for DRP computation with two modifications.  One problem is overflow of the initial value $\drpo{0}{a}$ for high values of the parameter $\beta$. When $\beta$ is integer, we can compute $\drpo{0}{a}$ by recurrence
\begin{equation}
\label{Racah_recurrence_R0}
\begin{array}{l}
\displaystyle
F(0)=\frac{\alpha+1}{(a+b)(\alpha+b-a)}\\
\displaystyle
F(k)=\frac{(\alpha+k+1)(2a-k+1)}{(a+b-k)(b-a+\alpha+k)}F(k-1),\hspace{4em} k=1,2,\ldots,\beta\\
\drpo{0}{a}=\sqrt{F(\beta)(2a+1)}.
\end{array}
\end{equation}
The factor $(\alpha+k)/(\alpha+k)$ from the original paper can be omitted.
The other values are obtained by the recurrence relation over $n$ as in Eq.~(\ref{Racah_recurrence_n}). It is called Algorithm~1.

Another algorithm is based on the recurrence over $s$. It begins by the same way, computation of $\drpo{0}{a}$ by Eq.~(\ref{Racah_recurrence_R0}). The initial values of higher degrees are
\begin{equation}
\label{Racah_recurrence_Rn0}
\begin{array}{l}
\drpo{n}{a}=\frac{(a-b+n)(\beta+n)(a+b+\alpha+n)}{n}\sqrt{D}\drpo{n-1}{a},\\
D=\frac{n(\alpha+\beta+2n+1)(\alpha+\beta+n)}{(\alpha+n)(\beta+n)(b-a+\alpha+\beta+n)(a+b+\alpha+n)(\alpha+\beta+2n-1)(a+b-\beta-n)(b-a-n)}.
\end{array}
\end{equation}
In the paper, there is incorrect factor $(\alpha+b-\beta-n)$ instead of $(a+b-\beta-n)$ in the denominator of $D$. The rest of the initial values is computed as
\begin{equation}
\label{Racah_recurrence_Rn1_wrong}
\drpo{n}{a+1}\!=\!E\sqrt{\frac{\rho(a+1)}{\rho(a)}\cdot\frac{2a+3}{2a+1}}\drpo{n}{a},\\
\end{equation}
where
\begin{equation}
\label{Racah_recurrence_Rn1_E_wrong}
E=\left(1+\frac{\scriptstyle{2n(\alpha+\beta+n+1)(a+1)}}{\scriptstyle{(a-b+1)(\beta+1)(a+b+\alpha+1)}}\right)
\end{equation}
and
\begin{equation}
\label{Racah_recurrence_Rn1_D_wrong}
\frac{\scriptstyle{\rho(a+1)}}{\scriptstyle{\rho(a)}}=\frac{\scriptstyle{(2a+1)(\beta+1)(b+\alpha+a+1)(b-a+1)}}{\scriptstyle{(b+\alpha-a-1)(2a-\beta+1)(a+b+1)}}.
\end{equation}
It is incorrect, the correct version is
\begin{equation}
\label{Racah_recurrence_Rn1_E}
\begin{array}{rcl}
E&=&\left(1-\frac{2\lambda(a+1)}{\tau(a)}\right)=\\
 &=&\left(1+\frac{2n(\alpha+\beta+n+1)(a+1)}
{(\alpha+1)(\beta+1)+a(a+1)(\alpha+\beta+2)
-a(\alpha+1)(a-\beta)-b(\beta+1)(b+\alpha)}
\right),
\end{array}
\end{equation}
i.e. the denominator is completely incorrect
and
\begin{equation}
\label{Racah_recurrence_Rn1_D}
\frac{\scriptstyle{\rho(a+1)}}{\scriptstyle{\rho(a)}}=\frac{\scriptstyle{(2a+1)(\beta+1)(b+\alpha+a+1)(b-a-1)}}{\scriptstyle{(b+\alpha-a-1)(2a-\beta+1)(a+b+1)}},
\end{equation}
i.e. in the numerator, there should be $(b-a-1)$ instead of $(b-a+1)$.

Finally, Daoui et al. use the stabilizing condition. When $\drpo{n}{s}$ is computed by Eq.~(\ref{Racah_recurrence_s}),
the new value is tested. When
\begin{equation}
\label{Racah_recurrence_Rn1_cond}
n>\frac{N}{6} \, \wedge \, \abs{\drpo{n}{s}}<10^{-6} \, \wedge \, \abs{\drpo{n}{s}}>\abs{\drpo{n}{s-1}},
\end{equation}
the value of $\drpo{n}{s}$ is substituted by zero. The symbol $\wedge$ means the logical and. It erases senselessly high values distorted by propagated error.
It is called Algorithm~3 in the paper. We will use it, after the error corrections, as the reference algorithm.

\subsubsection{Gram-Schmidt Orthogonalization}
Gram-Schmidt orthogonalization process (GSOP) is a way, how to change a set of functions to another set of orthogonal functions. It can be used for derivation of completely new orthogonal polynomials, e.g. GSOP applied on a set $\{1,x,x^2,\ldots\}$ in the interval $\langle -1,1 \rangle$ gives Legendre polynomials, see e.g.\cite{ThomasGSOP}. We can use GSOP also for increasing precision of orthogonal polynomials computed by another method. Here we have computed $\drpo{n}{s}$, but we are not sure, if it is sufficiently precise.  We can compute correction
\begin{equation}
\label{GSOP_cor}
\mathbfcal{T}(s)=\sum\limits_{k=0}^{n-1}\drpo{k}{s}\left(\sum\limits_{i=a}^{a+N-1}\drpo{n}{i}\drpo{k}{i}\right), \hspace{2em}s=a,a+1,\ldots,a+N-1.
\end{equation}
This correction is then subtracted from the original value
\begin{equation}
\label{GSOP_sub}
\check{\mathbfcal{R}}_{n}\!\left(s\right)=\drpo{n}{s}-\mathbfcal{T}(s), \hspace{2em}s=a,a+1,\ldots,a+N-1.
\end{equation}
Then we must correct also the norm
\begin{equation}
\label{GSOP_norm}
\tilde{\mathbfcal{R}}_{n}\!\left(s\right)=\check{\mathbfcal{R}}_{n}\!\left(s\right)/\left(\sqrt{\sum\limits_{i=a}^{a+N-1}{\check{\mathbfcal{R}}_{n}\!\left(s\right)}^2}+\varepsilon\right), \hspace{2em}s=a,a+1,\ldots,a+N-1,
\end{equation}
where $\varepsilon$ is some small value preventing division by zero. In Matlab $\varepsilon=2.2204\cdot 10^{-16}$.
$\tilde{\mathbfcal{R}}_{n}\!\left(s\right)$ is now version of $\drpo{n}{s}$ with increased precision.

GSOP works well, its main disadvantage is the high computing complexity $\mathcal{O}(N^3)$ (if we compute all the degrees up to $n=N-1$), while the computing complexity of all other algorithms mentioned in this paper is $\mathcal{O}(N^2)$.
It is big limitation of this method, the computing time may not be acceptable for very high $N$.

\subsection{The definition of discrete Racah moments (DRM)}
DRMs represent the projection of a signal (speech or images) on the basis of DRP.
The computation of the DRMs ($ \phi_{nm} $) for a 2D signal, $ f(x,y) $, with a size of $N_1\times N_2$ is performed by

\begin{align}
	\label{Eq_2D}
	&\phi_{nm} = \sum\limits_{x=0}^{N_1-1} \sum\limits_{y=0}^{N_2-1} f(x,y) \drpnn{1}{n}{x;N_1} \drpnn{2}{m}{y;N_2}  \\
	&n = 0,1,\dots,N_1-1 \ \ \mathrm{and}\ \
	m = 0,1,\dots,N_2-1. \nonumber		
\end{align}

The reconstruction of the 2D signal (image) from the Racah domain (the space of Racah moments) into the spatial domain can be carried out by

\begin{align}
	\label{Eq_2D_Rec}
	&\hat{f}(x,y)=\sum\limits_{n=0}^{N_1-1}\sum\limits_{m=0}^{N_2-1} \phi_{nm} \drpnn{1}{n}{x;N_1} \drpnn{2}{n}{y;N_2}  \\
	&\hspace{2em}x = 0,1,\dots,N_1-1 \ \ \mathrm{and}\ \
	y = 0,1,\dots,N_2-1. \nonumber 
\end{align}

\section{The Proposed Methodology}
\label{Method}
\begin{figure}[ht]
	\centering
	\includegraphics[width=0.99\linewidth]{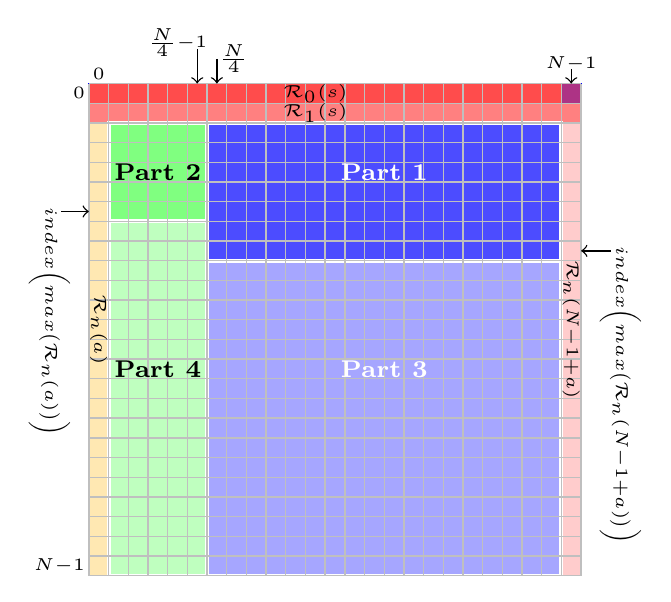}
	\caption{The matrix of DRPs. Note: the matrix is given by $n$ and $x$, where $x=s-a$.}
	\label{fig:DRP_plane}
\end{figure}

This section presents the proposed methodology for computing DRPs. We call it improved stabilization (ImSt).
\begin{align}
&\drpo{n}{s}= \sqrt{\frac {\Gamma\!\left( a{+}N{+}\alpha{-}s \right)\! \Gamma\! \left( a{+}N{+}\alpha{+}s{+}1 \right)  \left( \alpha{+}1{+}2n \right)\! \Gamma\! \left( N{-}n \right)\! \Gamma\! \left( 2a{+}N{-}n \right)  \left( 2s{+}1 \right) }{\Gamma \left(a{+}N{-}s \right) \Gamma \left( a{+}N{+}s{+}1 \right) \Gamma \left( N{+}\alpha{+}1{+}n \right) \Gamma \left( 2a{+}N{+}\alpha{+}n{+}1 \right) }} \times  \nonumber\\ 
& \hspace{2cm} \left( {-}N{+}1\right)_n \left( 2a{+}N{+}\alpha{+}1\right)_n {}_4F_3\Hyp{{-}n,a{-}s,a{+}s{+}1,\alpha{+}n{+}1}{1,2a{+}N{+}\alpha{+}1,{-}N{+}1}{1}
\end{align}
The DRP matrix is partitioned into four parts. They are shown in \figurename{ \ref{fig:DRP_plane}} as Part~1, Part~2, Part~3, and Part~4.
In the following subsections, the detailed steps are given. First of all, we must compute initial values.

\subsection{The First Initial Value}
The selection of the first initial value, specifically  its location and how it is computed, is considered crucial because all the other values of the the polynomial rely on that initial value.

The computation of the initial value in the existing algorithms limits the ability to compute the entire values of DRPs. For example, in \cite{daoui2022}, the formula for computation of the first initial value is as follows

\begin{align}
    &\drpo{0}{a}=\sqrt{F \times (2a+1)}\,, \nonumber \\
    &F=\frac{\Gamma(2a{+}1)\Gamma(\alpha{+}\beta{+}2)\Gamma(b{-}a{+}\alpha)\Gamma(a+b-\beta)}{\Gamma(a{+}b{+}1)\Gamma(\alpha{+}1)\Gamma(2a{+}1-\beta)\Gamma({-}a{+}b{+}\alpha{+}\beta{+}1)}\,. \label{Eq_R00}
\end{align}

\begin{figure*}[ht]
	\centering
	\subfloat[]{\includegraphics[width=1\columnwidth]{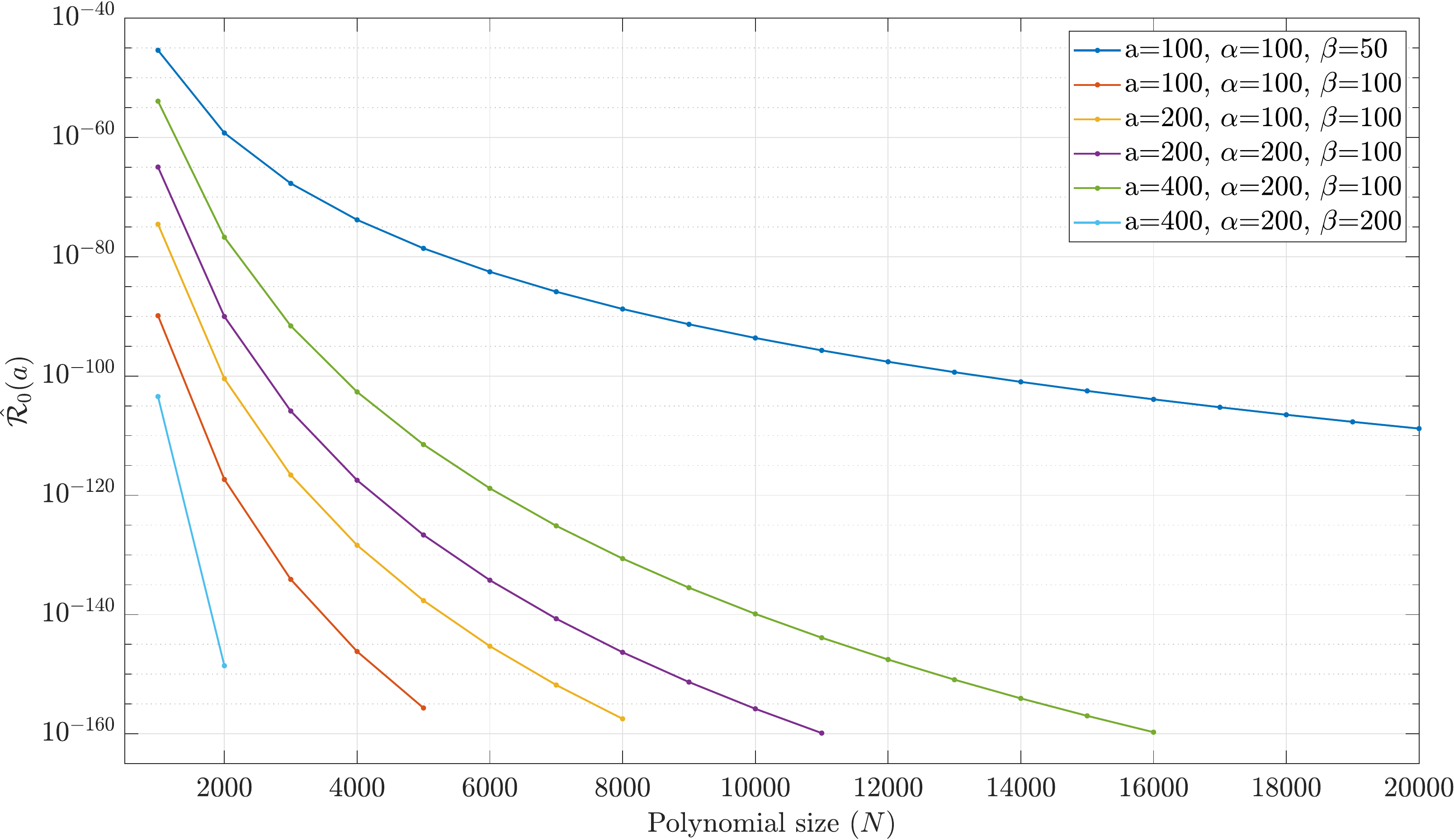}}
	\hfil
	\subfloat[]{\includegraphics[width=1\columnwidth]{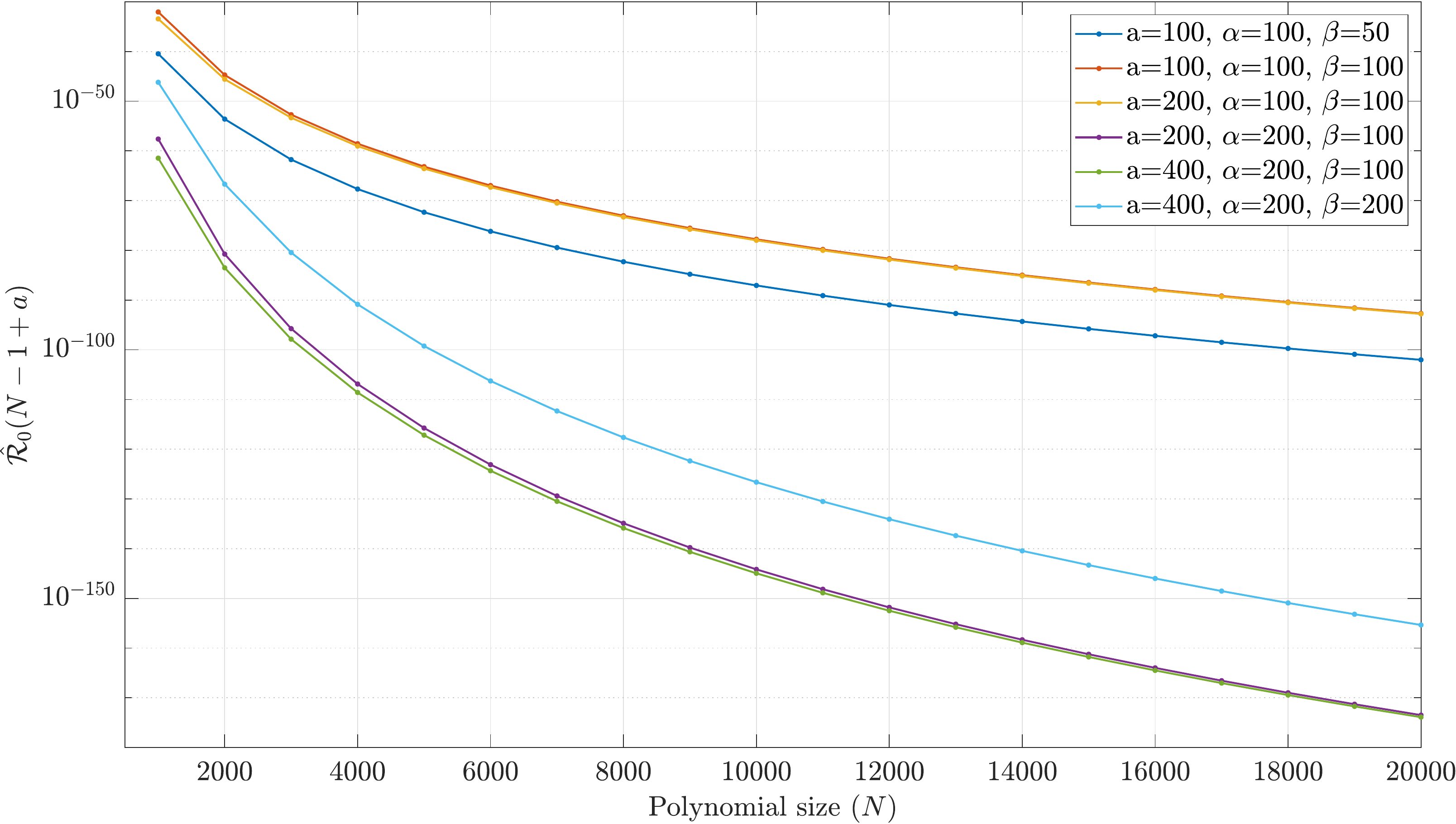}}
	\caption{The plot of the initial value (a) $\drpo{0}{a}$, and (b) $\drpo{0}{N-1+a}$.}
	\label{fig:initials}
\end{figure*}

This formula \eqref{Eq_R00} is uncomputable for a wide range of parameter values $a$, $\alpha$, and $\beta$ as shows in \figurename{ \ref{fig:initials}}a. Thus, in the proposed algorithm, we begin the computation at the last value of the first row, i.e. at $s=a+N-1$ as follows
\begin{align}
    &\drpo{0}{N-1+a}= \sqrt{\frac{\Gamma(\alpha{+}\beta{+}2)\Gamma(2a{+}N)\Gamma(\beta{+}N)\Gamma(2a{+}2N{+}\alpha)}{\Gamma(2a{+}2N-1)\Gamma(\beta{+}1)\Gamma(\alpha{+}\beta{+}N+1)\Gamma(2a{+}N{+}\alpha{+}1)}}, \label{Eq_R0N}
\end{align}
however, the Gamma function ($\Gamma(\cdot)$) make this equation uncomputable. To fix this issue, we rewrite Equation~\eqref{Eq_R0N} as

\begin{align}
    &Y=\psi(\alpha{+}\beta{+}2){+}\psi(2a{+}N){+}\psi(\beta{+}N){+}\psi(2a{+}2N{+}\alpha)-\nonumber \\
    &\hspace{2em}(\psi(2a{+}2N{-}1){+}\psi(\beta{+}1){+}\psi(\alpha{+}\beta{+}N+1){+}\psi(2a{+}N{+}\alpha{+}1)) \label{Eq_R0NF} \\
    &\drpo{0}{a+N-1}=\exp(Y/2)\nonumber
\end{align}
where $\psi(\cdot)$ represents the logarithmic gamma function: $\psi(x)=\log(\Gamma(x))$, $\log(\cdot)$ is natural logarithm. Using \eqref{Eq_R0NF} the first initial value is computable for a wide range of the DRP parameters as shown in \figurename{ \ref{fig:initials}}b.

\subsection{The Initial sets}
After computing the first initial value, the initial sets in the first two rows $\drpo{0}{s}$ and $\drpo{1}{s}$ are computed by the two-term recurrence relation. These initial sets will be used for computation of the remaining coefficients of DRPs (the coefficients in Parts 1, 2, 3, and 4).
The values of the coefficients $\drpo{0}{s}$ are calculated as follows
\begin{align}\label{Eq_R0s}
\drpo{0}{s}=&\sqrt{{\frac{\left( 2s{+}1\right)\left(a{-}{\beta}{+}s{+}1 \right)\left(b{+}s{+}1\right)\left(b{+}{\alpha}{-}s{-}1\right)\left(a{-}s{-}1 \right)}{\left(a{+}s{+}1\right)\left(b{+}{\alpha }{+}s{+}1\right)\left(a{-}{\beta}{-}s{-}1 \right)\left(2s{+}3\right)\left(b{-}s{-}1\right) }}}\  \drpo{0}{s{+}1}, \\
&s=a+N-2,a+N-3,\dots,a. \nonumber
\end{align}
After computation of the values $\drpo{0}{s}$, the values $\drpo{1}{s}$ are computed using the previously computed coefficients. The values of the coefficients of $\drpo{1}{s}$ are computed
\begin{align}
&\drpo{1}{s}= -\big(\left( (-a+b-1)\alpha+b^2-s^2-a-s-1 \right)\beta + (a^2-s^2+b-s-1)\alpha + a^2+b^2-2(s^2+s)-1 \big) \times \nonumber \\
& \hspace{1cm}\sqrt {{\frac {\alpha{+}\beta{+}3}{
(a{-}b{+}1)(a{+}b{-}\beta{-}1)(\alpha{+}1)(\beta{+}1)(a{-}b{-}\alpha{-}\beta{-}1)(a{+}b{+}\alpha{+}1)
}}}\times \drpo{1}{s{+}1}, s=a, a+1,\dots,a+N-1. \label{Eq_R1s}
\end{align}

\subsection{The Controlling Indices in the First and Last Columns}
\label{The_controlling_Indices}
To control the stability of the computation of the DRP coefficients, we present a controlling indices that are used to stabilize the computation of the coefficients. We first compute the coefficients $\drpo{n}{a}$ and $\drpo{n}{a+N-1}$ in the first and last columns. Then, the location, where the peak values occur, are found.
To compute the coefficients for $\drpo{n}{a}$, the two-term recurrence relation is used
\begin{align}
&\drpo{n+1}{a} = {-}\sqrt {{\frac { \left( N{-}n{-}1 \right)  \left( \alpha{+}\beta{+}2n{+}3  \right)  \left( \alpha{+}\beta{+}n{+}1 \right)  \left( \beta{+}n{{+}}1  \right)  \left( 2a{+}N{+}\alpha{+}n{+}1 \right) }{ \left( 2a{+}N{-}\beta {-}n{-}1 \right)  \left( \alpha{+}\beta{+}2\,n{+}1 \right)  \left( \alpha {+}n{+}1 \right)  \left( N{+}\alpha{+}\beta{+}n{+}1 \right)  \left( n{+}1  \right) }}}\times \drpo{n}{a}, \label{Eq_Rna} \\  &\hspace{4em}n=1,2,\dots,N-2. \nonumber
\end{align}
Also, we present a new two-term recurrence relation to compute the coefficients of $\drpo{n}{a+N-1}$ as follows
\begin{align}
&\drpo{n+1}{a+N-1} = \sqrt {{\frac { \left( N{-}n{-}1 \right)  \left( \alpha{+}\beta{+}2n{+}3  \right)  \left( \alpha{+}\beta{+}n{+}1 \right)  \left( \alpha{+}n{+}1 \right)  \left( 2a{+}N{-}\beta{-}n{-}1 \right) }{ \left( 2a{+}N{+}\alpha{+}n{+}1  \right) \left( \alpha{+}\beta{+}2n{+}1 \right)  \left( \beta{+}n{+}1 \right)  \left( N{+}\alpha{+}\beta{+}n{+}1 \right)  \left( n{+}1 \right)   }}}  \times \drpo{n}{a+N-1}, \label{Eq_RnN_1}\\ 
&\hspace{4em}     n=1,2,\dots,N-2. \nonumber 
\end{align}
The peak value at the last column $s=a+N-1$, i.e. the index
\begin{equation}
\label{indN}
ind_{N-1}=\operatorname*{arg\,max}\limits_{n=0,1,\ldots N-1} \drpo{n}{a+N-1}
\end{equation}
then creates the border between Part 1 and Part 3, while the border between Part 2 and Part 4 is the peak value at the first column $s=a$, i.e. the index
\begin{equation}
\label{ind0}
ind_{0}=\operatorname*{arg\,max}\limits_{n=0,1,\ldots N-1} \drpo{n}{a}.
\end{equation}

\subsection{The Controlling Index in the Last Row}
\label{The_controlling_Indices_last_row}
We would also need the index $N_s$ as the border between Part 1 and Part 2. The ideal value would be the peak value at the last row. We cannot compute it directly because of underflow for high $N$, then we can use substitutional value $N_s=\lfloor N/4+0.5\rfloor$ as is written in \figurename{ \ref{fig:DRP_plane}}. The symbol $\lfloor\cdot\rfloor$ is the fuction floor, the index $N_s$ is rounded to the nearest integer.

There is another possibility. Some values of $\drpo{n}{s}$ can underflow for high $N$, but the ratio of the adjacent values does not, therefore we can compute $N_s$ in logarithms. There is one complication, we need the logarithm of a sum $\log(a+b)$, but when $a$ and $b$ are similar, we can compute it as $\log((a/b+1)b)$ = $\log(a/b+1)+\log(b)$ = $\log(\exp(\log(a)-\log(b))+1)+\log(b)$.
The whole algorithm is then as follows. First, we compute logarithm of the first value
\begin{equation}
\label{Eq_l00}
\begin{array}{l}
L_{0,0}=\log(\drpo{0}{a})=(\psi(2a+2)+\psi(N+\alpha)+\psi(\alpha+\beta+2)+\psi(2a+N-\beta)-\\ \ \hspace{2cm}\psi(2a-\beta+1)-\psi(2a+N+1)-\psi(\alpha+1)-\psi(N+\alpha+\beta+1))/2.
\end{array}
\end{equation}
Again, $\psi(x)$ is the logarithmic gamma function. Then we compute values in the first column. We need not remember them, we need only the last value $L_{N,0}$.
\begin{equation}
\label{Eq_ln0}
\begin{array}{l}
L_{n,0} = L_{n-1,0} + (\log(N{-}n{-}1){+}\log(\alpha{+}\beta{+}2n{+}3){+}\log(\alpha{+}\beta{+}n{+}1){+}\log(\beta{+}n{+}1){+}\\ \hspace{2cm}\log(2a{+}N{+}\alpha{+}n{+}1){-}\log(2a{+}N{-}\beta{-}n{-}1){-}\log(\alpha{+}\beta{+}2n{+}1){-}\\\hspace{2cm}\log(\alpha{+}n{+}1){-}\log(N{+}\alpha{+}\beta{+}n{+}1){-}\log(n{+}1))/2\\
n=1,2,\ldots,N.
\end{array}
\end{equation}
The signum of the result must be computed separately
\begin{equation}
\label{Eq_sn0}
S_{N,0}=(-1)^{(N-1)}.
\end{equation}
The second value in the last row
\begin{equation}
\label{Eq_ln1}
\begin{array}{l}
L_{N,1}=\log(\lvert LE\rvert)+\log(LD)/2+L_{N,0} \\
S_{N,1}=\mathrm{sign}(LE)S_{N,0},
\end{array}
\end{equation}
where
\begin{equation}
\label{Eq_ln1ad}
\begin{array}{l}
\displaystyle
LE=1-\frac{2(N{-}1)(\alpha{+}\beta{+}n{+}1)(a{+}1)}{a(\alpha{+}1)(a{-}\beta){+}b(\beta{+}1)(b{+}\alpha){-}(\alpha{+}1)(\beta{+}1){-}a(a{+}1)(\alpha{+}\beta{+}2)}\\
\displaystyle
LD=\frac{(2a{+}1)(\beta{+}1)(a{+}b{+}\alpha{+}1)(b{-}a{-}1)(2a{+}3)}{(b{-}a{+}\alpha{-}1)(2a{-}\beta{+}1)(a{+}b{+}1)(2a{+}1)}.
\end{array}
\end{equation}
The factor $LA$ equals $E$ from Eq.~\eqref{Racah_recurrence_Rn1_E} and
$$
LD=\frac{\rho(a+1)}{\rho(a)}\cdot\frac{2a+3}{2a+1}
$$
from Eq.~\eqref{Racah_recurrence_Rn1_D}.

The last row is then computed by the recurrence
\begin{equation}
\label{Eq_lnx}
\begin{array}{l}
L_{N,x}=\log(\lvert LA\rvert)+\log(\lvert B1\rvert)+\log(\lvert B\rvert)/2+L_{N,x-2}\\
S_{N,x}=S_1\cdot S_3\\
x=2,3,\ldots\\
s=a+x,
\end{array}
\end{equation}
where
\begin{align}
A&=\frac{(a+s)(s-a+\beta)(b+\alpha+s)(b-s)(2s+1)}{(b+\alpha-s)(a-\beta+s)(s-a)(s+b)(2s-1)}\nonumber\\
B&=A\frac{(a+s-1)(s-a+\beta-1)(b+\alpha+s-1)(b-s+1)(2s-1)}{(b+\alpha-s+1)/(a-\beta+s-1)(s-a-1)(s+b-1)(2s-3)}\nonumber\\
A1&=\frac{(2s-1)(\sigma(s-1)+(s-1)\tau(s-1)-2\lambda s(s-1))}{(s-1)(\sigma(s-1)+(2s-1)\tau(s-1))}\nonumber\\
B1&=\frac{s\sigma(s-1)}{(s-1)(\sigma(s-1)+(2s-1)\tau(s-1))}\label{Eq_lnxw}\\
S_1&=\mathrm{sign}(B1)S_{N,x-2}\nonumber\\
S_2&=\mathrm{sign}(A1)S_{N,x-1}S_1\nonumber\\
LA&= S_2\exp(\log(\lvert A1\rvert)+\log(A)/2+L_{N,x-1}-\log(\lvert B1\rvert)-\log(B)/2-L_{N,x-1}-1\nonumber\\
S_{3}&=\mathrm{sign}(LA),\nonumber
\end{align}
and the functions $\sigma(s)$, $\tau(s)$, and $\lambda$ are given in Eq.~\eqref{Racah_recurrence_s_func}.
When we find the maximum, i.e. the point $x$, where $L_{N,x-1}>L_{N,x}$ $\wedge$ $L_{N,x-1}>L_{N,x-2}$, then
we have found the index $N_s=x-1=s-a-1$. It is better to stop the computation here, because if $\sigma(s)$ is higher than about $4\cdot 10^{15}$, the computation looses its accuracy.
In our tests it was always after the finding $N_s$.

It is also possible to compute the maximum from the end of the last row. The value $L_{0,N}$ equals $Y/2$, where $Y$ is from Eq.~\eqref{Eq_R0NF}.
Then we compute the values in the last column because of the last value
\begin{equation}
\label{Eq_lnN}
\begin{array}{l}
L_{n,N} = L_{n-1,N} +(\log(N-n-1)+\log(\alpha+\beta+2n+3)+\log(\alpha+\beta+n+1)+\\
\hspace{2cm}\log(\alpha+n+1)+\log(2a+N-\beta-n-1)-\log(\alpha+\beta+2n+1)-\\
\hspace{2cm}\log(\beta+n+1)-\log(N+\alpha+\beta+n+1)-\log(n+1)-\\
\hspace{2cm}\log(2a+N+\alpha+n+1))/2,
\hspace{10em}n=1,2,\ldots,N.
\end{array}
\end{equation}
The last but one value in the last row
\begin{equation}
\label{Eq_ln1mx}
\begin{array}{l}
L_{N,N-1}=L_{N,N}+\log(\lvert LF\rvert)-\log(LG)/2 \\
S_{N,N-1}=\mathrm{sign}(LF),
\end{array}
\end{equation}
where
\begin{equation}
\label{Eq_ln1admx}
\begin{array}{l}
\displaystyle
LF=\frac{(2s-1)(\sigma(s-1)+(s-1)\tau(s-1)-2\lambda_ns(s-1))}{s\sigma(s-1)}\\
\displaystyle
LG=\frac{(a+s-1)(s-a+\beta-1)(b+\alpha+s-1)(b-s+1)(2s-1)}{(b+\alpha-s+1)(a-\beta+s-1)(s-a-1)(s+b-1)(2s-3)}\\
s=a+N.\\
\end{array}
\end{equation}
We can invert the recurrence for the direct computation
\begin{equation}
\label{Eq_lnx1}
\begin{array}{l}
L_{N,x-2}\!=\!\log(\lvert LA\rvert)-\log(\lvert B1\rvert)-\log(\lvert B\rvert)/2+L_{N,x-1}\\
S_{N,x-2}\!=\!S_1\cdot S_3\cdot\mathrm{sign}(B1)\\
x=N-2,N-3,\ldots\\
\end{array}
\end{equation}
where $A$, $B$, $A1$, $B1$ is the same as in Eq.~\eqref{Eq_lnx} and
\begin{equation}
\label{Eq_lnxi}
\begin{array}{l}
S_1=S_{N,x}\\
S_2=\mathrm{sign}(A1)S_{N,x-1}S_1\\
LA=S_2\exp(\log(\lvert A1 \rvert)+\log(A)/2+L_{N,x-1}-L_{N,x})-1\\
S_{3}=\mathrm{sign}(LA).
\end{array}
\end{equation}
The peak value $N_n$ is then the first value $x$, where
$L_{N,x-1}>L_{N,x}$ $\wedge$ $L_{N,x-1}>L_{N,x-2}$. Then $N_n=x-1=s-a-1$. Again, we should stop the computation here. If $N_n=N_s$, it is good indication that we have the correct value.

\subsection{The coefficients for Parts 1 and 2}
The coefficients in Parts 1 and 2 are computed using the three-term recurrence algorithm in the $n$-direction as follows
\begin{equation}\label{Eq_nr}
    \drpo{n}{s}=\Theta_1 \, \drpo{n-1}{s} + \Theta_2 \, \drpo{n-2}{s},
\end{equation}
where
\begin{align}
    \Theta_1 = \frac{\Theta_{11}}{\Theta_0}\sqrt{\Theta_{12}},\hspace{2em}
    \Theta_2 = \frac{\Theta_{21}}{\Theta_0}\sqrt{\Theta_{12}\Theta_{22}}
\end{align}
and
\begin{align}
    &\Theta_0= \frac{n \left(\alpha+\beta+n \right)}{\left(\alpha+\beta+2n-1\right)  \left(\alpha+\beta+2n \right)} \\
    &\Theta_{11}=  s(s{+}1) {-} \frac{1}{4}\left({a}^{2}{+}{b}^{2}{+} \left( a{-}\beta \right) ^{2}{+}\left( b{+}\alpha \right) ^{2}{-}2\right) {+}\\&\hspace{2cm} \frac{1}{8}
    \left( \left(\alpha{+}\beta+2n{-}2 \right)\left(\alpha{+}\beta{+}2n\right) \right)- \frac{1}{2} \left(
    \frac {\left({\beta}^{2}{-}{\alpha}^{2}\right)\left( \left( b{+}\alpha /2 \right)^{2}{-} \left( a{-}\beta /2 \right)^{2} \right)}
    {\left( \alpha{+}\beta{+}2n{-}2 \right)  \left( \alpha{+}\beta{+}2n \right)} \right) \nonumber
\end{align}
\begin{align}
    \Theta_{21}=& {-}{\frac { \left( \alpha{+}n{-}1 \right)  \left( \beta{+}n{-}1 \right) }{ \left( \alpha+\beta{+}2n{-}2 \right)  \left( \alpha+\beta{+}2 n{-}1 \right) }}\! \left( \left( a{+}b{+}\frac{\alpha{-}\beta}{2}\right) ^{2}\!\!{-}\! \left( n{-}1{+}\frac{\alpha+\beta}{2}  \right)^{2} \right) \!\times \nonumber \\ & \hspace{2cm} \left( \left( b{-}a{+}\frac{\alpha+\beta}{2}\right) ^{2}{-} \left( n{-}1{+}\frac{\alpha{+}\beta}{2} \right) ^{2} \right)
\end{align}
\begin{align}
    \Theta_{12}=& \frac {n \left( \alpha+\beta+n  \right) \left( \alpha+\beta+2n+1 \right) }{ \left( \alpha+n \right)  \left( \beta+n \right) \left( \alpha+\beta+2n-1 \right) \left( a-b-\alpha-\beta-n \right) \left( a-b+n \right) }\times\frac{1}{ \left( a+b+\alpha+n \right) \left( a+b-\beta-n \right) }
\end{align}
\begin{align}
    \Theta_{22}=& \frac { \left( n-1 \right)  \left( \alpha+\beta+n-1 \right)  \left( \alpha+\beta+2n-1 \right) }{ \left( \alpha{+}n{-}1 \right) \left( \beta{+}n{-}1 \right) \left( \alpha{+}\beta{+}2n{-}3 \right) \left( a{-}b{-}\alpha{-}\beta{-}n{+}1 \right) \left( a{-}b{+}n{-}1 \right)}\times\frac{1}{ \left( a{+}b{+}\alpha{+} n{-}1 \right) \left( a{+}b{-}\beta{-}n{+}1 \right) }.
\end{align}
The border between Part 1 and Part 3 is the index $ind_{N-1}$, see Eq.~\eqref{indN}, the recurrence algorithm is applied for $s=a+N_s,a+N_s+1,\dots,a+N-2$ and $n=2,3,\dots,ind_{N-1}$, while the border between Part 2 and Part 4 is the index $ind_{0}$, see Eq.~\eqref{ind0}, the recurrence algorithm is applied for $s=a+1,a+2,\dots,a+N_s-1$ and $n=2,3,\dots,ind_{0}$.

\subsection{The coefficients for Parts 3 and 4}
The coefficients in Parts 3 and 4 are computed using the same three-term recurrence algorithm in the $n$-direction as in \eqref{Eq_nr}.
After computation of each value, the following stabilizing condition is applied for each order $n$

\begin{equation}
    \drpo{n}{s}=0 \text{ if }
    \left\lvert\drpo{n}{s}\right\rvert{<}10^{-5} \wedge \left\lvert\drpo{n}{s}\right\rvert{>}\left\lvert\drpo{n-1}{s}\right\rvert.
\end{equation}
In Part 3, we add a codition, that there must exist $\drpo{i}{s}$ such that $\lvert\drpo{i}{s}\rvert<10^{-5}$ for some $i<n$.

The recurrence algorithm for Part 3 is applied in the range $s=a+N_s,a+N_s+1,\dots,a+N-2$ and $n=ind_{N-1}+1,ind_{N-1}+2,\dots,N-1$; while for Part 4, the recurrence algorithm is carried out in the range $s=a+1,a+2,\dots,a+N_s-1$ and $n=ind_{0},ind_{0}+1,\dots,N-1$.

\subsection{Special case of Racah Polynomials}
\begin{figure}[!htt]
	\centering
	\includegraphics[width=0.99\linewidth]{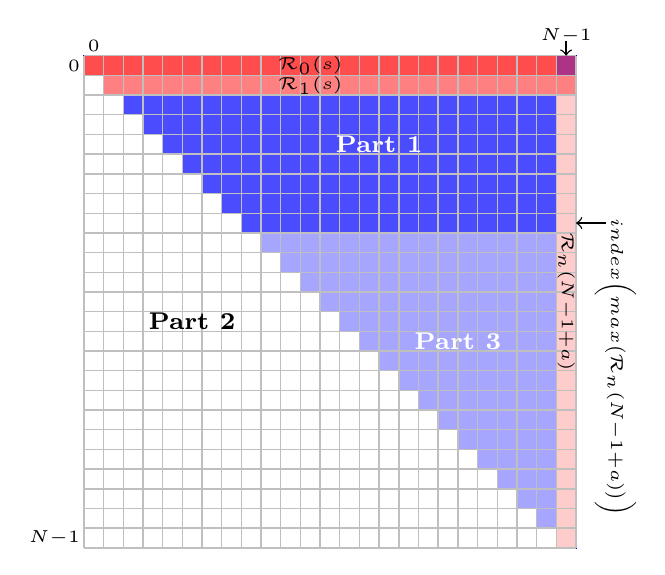}
	\caption{The matrix of DRPs for the case $a=\alpha=\beta=0$. }
	\label{fig:DRP_plane2}
\end{figure}

In this section, a special case of DRPs is presented.
The parameter $\beta$ affects on the energy compaction as its value becomes larger than 0.
So, the case $\drps{n}{s}$, where $a=\alpha=\beta=0$, has special significance.
In this case, the $\drps{n}{s}$ is given as follows
\begin{align}
    \drps{n}{s}=& \frac{(b{+}1)_n(1)_n({-}b{+}1)_n}{n!} \times
	{}_4F_3\Hyp{{-}n,{-}s,s{+}1,n{+}1}{1,b{+}1,{-}b{+}1}{1} \times\sqrt{\frac{\frac{\Gamma(s+1)\Gamma(b+s+1)\Gamma(b-s)\Gamma(s+1)}{\Gamma(b+s+1)\Gamma(b-s)\Gamma(s+1)\Gamma(s+1)}}{\frac{\Gamma(n+1)\Gamma(n+1)\Gamma(b+n+1)\Gamma(b+n+1)}{(2n+1)\Gamma(n+1)\Gamma(b-n)\Gamma(n+1)\Gamma(b-n)}}(2s{+}1)} \nonumber\\
	=&\frac{(b{+}1)_n(1)_n({-}b{+}1)_n}{n!} \times
	{}_4F_3\Hyp{{-}n,{-}s,s{+}1,n{+}1}{1,b{+}1,{-}b{+}1}{1} \times \sqrt{\frac{(2n+1)\Gamma(b-n)\Gamma(b-n)}{\Gamma(b+n+1)\Gamma(b+n+1)}(2s{+}1)} \nonumber\\
	=&\frac{(b{+}1)_n(1)_n({-}b{+}1)_n}{n!} \times	{}_4F_3\Hyp{{-}n,{-}s,s{+}1,n{+}1}{1,b{+}1,{-}b{+}1}{1} \times \frac{\Gamma(b-n)}{\Gamma(b+n+1)} \sqrt{(2n+1)(2s{+}1)} \nonumber\\
	=&\frac{\Gamma(b+n+1)n!\Gamma(-b+1+n)}{n!\Gamma(b+1)\Gamma(-b+1)} \times {}_4F_3\Hyp{{-}n,{-}s,s{+}1,n{+}1}{1,b{+}1,{-}b{+}1}{1} \times \frac{\Gamma(b-n)}{\Gamma(b+n+1)} \sqrt{(2n+1)(2s{+}1)} \nonumber\\
	=&\frac{\Gamma(b-n)\Gamma(-b+1+n)}{\Gamma(b+1)\Gamma(-b+1)} \times {}_4F_3\Hyp{{-}n,{-}s,s{+}1,n{+}1}{1,b{+}1,{-}b{+}1}{1} \times \sqrt{(2n+1)(2s{+}1)}.\label{Eq_aab}
\end{align}
Roman \cite{Roman1992} shows the property of factorial
\begin{equation}\label{Eq_factorial}
    c!(-c-1)! = (-1)^{c+(c<0)},
\end{equation}
where
\begin{equation}
    (c<0)=\left\{\begin{matrix}
1 & \text{ if } c<0\  \\
0 & \text{ if } c\geq 0. \\
\end{matrix}\right.
\end{equation}
It is well known that $c! = \Gamma(c+1)$; thus \eqref{Eq_factorial} can be written by this way
\begin{equation}\label{Eq_Gamma}
    \Gamma(c+1)\Gamma(-c) = (-1)^{c+(c<0)}.
\end{equation}
Using \eqref{Eq_Gamma}, the term $\Gamma(b-n)\Gamma(-b+1+n)$ from \eqref{Eq_aab} can be expressed
\begin{equation}
\label{Eq_num1}
\Gamma(b{-}n)\Gamma({-}b{+}1{+}n){=}\Gamma(b{-}n)\Gamma({-}(b{-}n){+}1){=}({-}1)^{{-}(b{-}n){+}1}{=}{-}({-}1)^{{-}b}({-}1)^n.
\end{equation}
Also, the term $\Gamma(b+1)\Gamma(-b+1)$ from \eqref{Eq_aab} can be expressed
\begin{equation}
\label{Eq_den1}
\Gamma(b+1)\Gamma(-b+1)= \Gamma(b+1)\Gamma(-b)(-b)= (-b)(-1)^{b+0}= -b(-1)^{-b}.
\end{equation}
From \eqref{Eq_num1} and \eqref{Eq_den1}, \eqref{Eq_aab} can be expressed
\begin{align}
	&\drps{n}{s}=\frac{-(-1)^{-b}(-1)^n}{-b(-1)^{-b}}  \sqrt{(2n{+}1)(2s{+}1)}\  {}_4F_3\Hyp{{-}n,{-}s,s{+}1,n{+}1}{1,b{+}1,{-}b{+}1}{1} \nonumber \\
	&=\frac{(-1)^n \sqrt{(2n+1)(2s{+}1)}}{b}\ \  {}_4F_3\Hyp{{-}n,{-}s,s{+}1,n{+}1}{1,b{+}1,{-}b{+}1}{1}. \label{Eq_aab1}
\end{align}
For \eqref{Eq_aab1}, replacing $n$ by $s$, we obtain
\begin{equation}
\label{Eq_aab2}
\drps{s}{n}=\frac{(-1)^s \sqrt{(2s+1)(2n{+}1)}}{b} \ \ {}_4F_3\Hyp{{-}s,{-}n,n{+}1,s{+}1}{1,b{+}1,{-}b{+}1}{1}.
\end{equation}
By comparing \eqref{Eq_aab1} with \eqref{Eq_aab2}, we obtain the following symmetry relation
\begin{equation}\label{Eq_sym}
    \drps{s}{n}=(-1)^{(s-n)}\drps{n}{s}.
\end{equation}
Thus, from \eqref{Eq_sym}, we can compute the coefficients for 50\% and the rest of the coefficients using the symmetry relation. In other words, the coefficients are computed in the range $n=0,1,\dots,N-1$ and $s=n,n+1,\dots,N-1$  (Parts 1 and 3). The rest of the coefficients are computed using the symmetry relation (Part 2) as shown in \figurename~\ref{fig:DRP_plane2}.
The Eq.~\eqref{Eq_nr} becomes
\begin{equation}\label{Eq_nr0}
    \drpo{n}{s}=\Theta_{10} \, \drpo{n-1}{s} + \Theta_{20} \, \drpo{n-2}{s},
\end{equation}
where
\begin{align}
    &\Theta_{10}=  \frac{\left(2s(s{+}1) {+} n\left( n{-}1\right){-}{N}^{2}{+}1\right)\sqrt{4n^2-1}}{n(N{-}n)(N{+}n)}\\
    &\Theta_{20}= {-}\frac { \left( n{-}1 \right) \left( N{-}n{+}1 \right) \left( N{+}n{-}1 \right) }{ n \left( N{-}n \right) \left( N{+}n \right)}\sqrt{\frac { 2n{+}1 }{ 2n{-}3 }}.
\end{align}

\subsection{Implementation of the proposed algorithm}
In this section the pseudo code is presented. The pseudo code of the proposed algorithm for the general case is presented in Algorithm \ref{algorithm1}. In addition, the pseudo code for the special case ($a=\alpha=\beta=0$) is given in Algorithm \ref{algorithm2}.
\begin{algorithm}
	\hspace*{\algorithmicindent} \textbf{Input: $Ord, a, b, \alpha, \beta$} \\
	\hspace*{\algorithmicindent} \hspace*{\algorithmicindent} \hspace*{\algorithmicindent} \text{$Ord$ is the maximum degree of DRP,  $Ord<b-a$.} \\
	\hspace*{\algorithmicindent} \hspace*{\algorithmicindent} \hspace*{\algorithmicindent} \text{$a, b, \alpha, \beta$ represents the parameter of DRP.} \\
	\hspace*{\algorithmicindent} \textbf{Output: $\drpo{n}{s}$}
	\begin{algorithmic}[1]
		\caption{Computation of the DRP coefficients using the proposed algorithm.} \label{algorithm1}
		
		\State $N \gets b-a$ \Comment{$N$ represents the size of DRP}
        \State $\Theta=10^{-5}$ \Comment{Threshold for stabilizing condition}

		\State Compute $\drpo{0}{N-1+a}$ using \eqref{Eq_R0NF}
		\For {$s=a+N-2 : a$}
		\State Compute $\drpo{0}{s}$ using \eqref{Eq_R0s}
		\EndFor
		\For {$s=a:a+N-1$}
		\State Compute $\drpo{1}{s}$ using \eqref{Eq_R1s}
		\EndFor
		\For {$n=1:Ord-1$}
		\State Compute $\drpo{n}{a}$ using \eqref{Eq_Rna}
		\State Compute $\drpo{n}{N-1+a}$ using \eqref{Eq_RnN_1}
		\EndFor

		\State $N_s$ from Sec.~\ref{The_controlling_Indices_last_row} or $N_s \gets \lfloor\frac{N}{4}+0.5\rfloor$ \Comment{$\lfloor x\rfloor$ is integer part of $x$}
		
		\State $ind_0 \gets\!\! \operatorname*{arg\,max}\limits_{n=0,1,\ldots N-1}\drpo{n}{a}$ \Comment Index of maximum $\drpo{n}{a}$
		
		\State $ind_{N-1} \gets\!\! \operatorname*{arg\,max}\limits_{n=0,1,\ldots N-1}\drpo{n}{a\!+\!N\!-\!1}$ \Comment Index of maximum $\drpo{n}{a\!+\!N\!-\!1}$
		
		\For {$s=a+N_s:a+N-1$} \Comment{Part 1}
		\For {$n=2:ind_{N-1}-1$}
		\State Compute $\drpo{n}{s}$ using \eqref{Eq_nr}
		\EndFor
		\EndFor
		
		\For {$s=a:a+N_s-1$} \Comment{Part 2}
		\For {$n=2:ind_0$}
		\State Compute $\drpo{n}{s}$ using \eqref{Eq_nr}
		\EndFor
		\EndFor

		\For {$s=a+N_s:a+N-1$} \Comment{Part 3}
		\For {$n=ind_{N-1}:Ord$}
		\State Compute $\drpo{n}{s}$ using \eqref{Eq_nr}
		\If{$\abs{\drpo{n}{s}}<\Theta \wedge \abs{\drpo{n}{s}}>\abs{\drpo{n-1}{s}} \wedge \exists \abs{\drpo{i}{s}}>\Theta,\ i<n$}
		\State $\drpo{n}{s}=0$
		\State Exit inner loop
		\EndIf
		\EndFor
		\EndFor

        \newcounter{lastline}
        \setcounter{lastline}{\value{ALG@line}}
	\end{algorithmic}
\end{algorithm}
\begin{algorithm}
    \begin{algorithmic}[1]
        \setcounter{ALG@line}{\value{lastline}}
		
		\For {$s=a+1:a+N_s-1$} \Comment{Part 4}
		\For {$n=ind_0+1:Ord$}
		\State Compute $\drpo{n}{s}$ using \eqref{Eq_nr}
		\If{$\abs{\drpo{n}{s}}<\Theta \wedge \abs{\drpo{n}{s}}>\abs{\drpo{n-1}{s}}$}
		\State $\drpo{n}{s}=0$
		\State Exit inner loop
		\EndIf
		\EndFor
		\EndFor		
	\end{algorithmic}
\end{algorithm}

\begin{figure}[ht]
\centering
\includegraphics[width=1\columnwidth,clip]{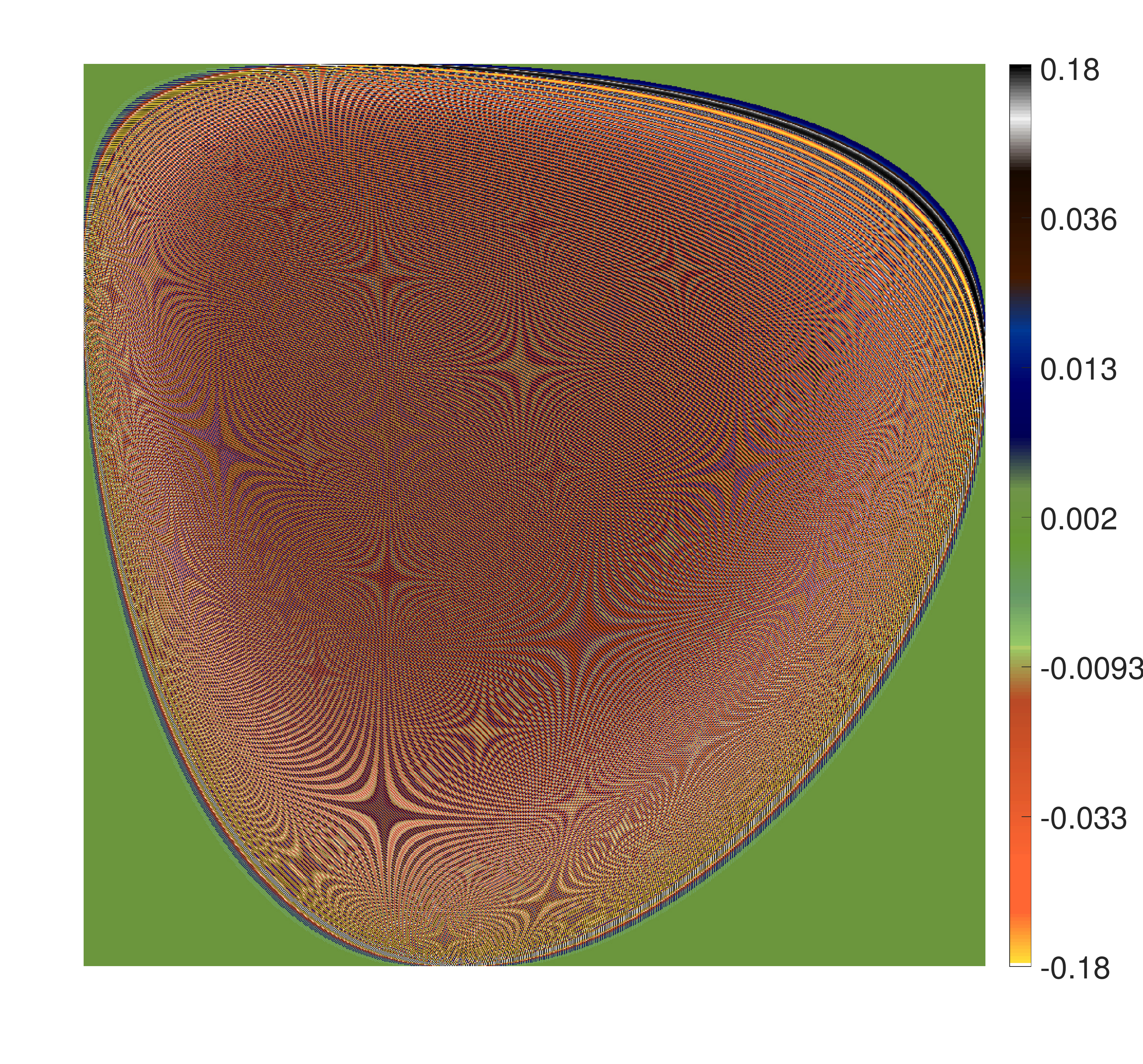}
\caption{The Racah polynomials for $a=800$, $b=1800$, $\alpha=400$, and $\beta=100$.}
\label{fig:values}
\end{figure}

\begin{algorithm}
	\hspace*{\algorithmicindent} \textbf{Input: $N, Ord$} \\
	\hspace*{\algorithmicindent} \hspace*{\algorithmicindent} \hspace*{\algorithmicindent} \text{$N$ represents the size of the DRP,} \\
	\hspace*{\algorithmicindent} \hspace*{\algorithmicindent} \hspace*{\algorithmicindent} \text{$Ord$ is the maximum degree of the DRP, $Ord<N$.} \\
	\hspace*{\algorithmicindent} \textbf{Output: $\drpo{n}{s}$}
	\begin{algorithmic}[1]
		\caption{Computation of the DRP coefficients using the proposed algorithm for the special case $a=\alpha=\beta=0$.} \label{algorithm2}
		
        \State $\Theta=10^{-5}$ \Comment{Threshold for stabilizing condition}
		\State $\drpo{0}{N-1} \gets \sqrt{2N-1}/N$
		\For {$s=N-2 : 0$}
		\State $\drpo{0}{s} \gets \sqrt{(2s+1)/(2s+3)}\times \drpo{0}{s+1}$
		\EndFor
		\For {$s=1 : Ord$}
		\State $\drpo{s}{0} \gets (-1)^{s}\drpo{0}{s}$
		\EndFor
		
		\For {$s=0:N-1$}
		\State $\drpo{1}{s} \gets -(N^2-2s^2-s2-1)\sqrt{3}/(N^2-1) \times \drpo{0}{s}$
		\EndFor
		\For {$s=2:Ord$}
		\State $\drpo{s}{1} \gets (-1)^{s-1}\drpo{1}{s}$
		\EndFor
		
		\For {$n=1:Ord-1$}
		\State $\drpo{n{+}1}{N{-}1} \gets  (N{-}n{-}1) \sqrt{2n{+}3}/(N{+}n{+}1)/\sqrt{2n{+}1}\times \drpo{n}{N{-}1}$
		\EndFor

		\State $ind_{N-1} \gets \!\!\operatorname*{arg\,max}\limits_{n=0,1,\ldots N-1}\drpo{n}{N\!-\!1}$ \Comment position of maximum in $\drpo{n}{N\!-\!1}$

		\For {$n=2:ind_{N-1}-1$} \Comment{Part 1}
		\For {$s=n:N-1$}
		\State Compute $\drpo{n}{s}$ using \eqref{Eq_nr0}
		\EndFor
		\EndFor
		
		\For {$n=ind_{N-1}:Ord$} \Comment{Part 3}
		\For {$s=n:N-1$}		
		\State Compute $\drpo{n}{s}$ using \eqref{Eq_nr0}
		\If{$\abs{\drpo{n}{s}}<\Theta  \wedge  \abs{\drpo{n}{s}}>\abs{\drpo{n-1}{s}}$}
		\State $\drpo{n}{s}=0$
		\State Exit inner loop
		\EndIf
		\EndFor
		\EndFor

		\For {$s=3:Ord$} \Comment{Part 2}
		\For {$n=2:s-1$}
		\State $\drpo{s}{n} \gets (-1)^{s-n}\drpo{n}{s}$
		\EndFor
		\EndFor
		
	\end{algorithmic}
\end{algorithm}
The values of the Racah polynomials for $a=800$, $b=1800$, $\alpha=400$, and $\beta=100$ (i.e. N=1000) in artificial colors are in \figurename{ \ref{fig:values}}.

\section{Experimental Analysis}
\label{Exper}
This section evaluates the proposed algorithm for DRP and compares it with the existing algorithms.
Three evaluation procedures are carried out to check the performance of the proposed algorithm which are: maximum size generated, computational cost, and signal reconstruction. The experiments were carried out using MATLAB version 2019b on the computer with the processor Intel(R) Core(TM) i9-7940X CPU with frequency 3.10GHz, memory	32,0 GB and with 64-bit Windows 10 Pro.

\subsection{Maximum Degree}
We searched the maximum signal size $N$, where the orthogonality error $E$ is less than 0.001. We changed the parameter values $a$, $\alpha$ and $\beta$ as ratios of $N$. It has an advantage, that the pattern of non-zero values looks similar and is not moved.
The orthogonality error is defined
\begin{align}
\label{Eq_oe}
E=\max\limits_{n,m=0,1,\ldots,N-1}\left\lvert\sum\limits_{s=a}^{b-1}\drpo{n}{s} \drpo{m}{s} - \delta_{nm}\right\rvert\ .
\end{align}
The results are in Tab.~\ref{Tab:Racah}.
\begin{table}[!httbp]
\begin{center}
\caption{Maximum sizes $N$ of the Racah polynomials reachable by various algorithms. Usually, the limit is the algorithm precision, i.e. the orthogonality error $E\leq 10^{-3}$, $\dag$ the limit is computing time $\leq$1 hour, $\ddag$ the limit is computer memory 32 GB.}
\label{Tab:Racah}
\vspace{2ex}
\begin{tabular}{|l|cccc|}
\hline
 & $a=0$ & $a=\lceil N/10000+0.5\rceil$ & $a=\lfloor N/4+0.5\rfloor$  & $a=\lfloor N/2+0.5\rfloor$ \\
 & $\alpha=0$ & $\alpha=N/10000$ & $\alpha=\lfloor N/8+0.5\rfloor$  & $\alpha=\lfloor N/2+0.5\rfloor$ \\
 & $\beta=0$ & $\beta=N/10000$ & $\beta=\lfloor N/16+0.5\rfloor$  & $\beta=\lfloor N/4+0.5\rfloor$ \\
\hline
 Zhu $n$ & 23 & 25 & 37 & 32   \\
 Zhu $s$ & 21 & 26 & 35 & 32   \\
 Daoui   & 1165 & 4 & 65 & 53   \\
 GSOP    & $9649^{\dag}$ & $9834^{\dag}$ & 1075 & 504  \\
 ImSt    & $56000^{\ddag}$ & 25580 & 6770 & 4659 \\
\hline
\end{tabular}
\end{center}
\end{table}

In the first column, when $a=0$, $\alpha=0$ and $\beta=0$, our algorithm~\ref{algorithm2} is used, in the other cases, it is our algorithm~\ref{algorithm1}.
The limit $N=56000$ is not limit of our algorithm, it is the memory limit of our computer. We are not able to check the orthogonality error because of the ``Out of memory'' error.

Another problem is long computation of GSOP. In the case $a=\alpha=\beta=0$ and $N=56000$, the error of orthogonality $E$ was also under the threshold 0.001, but the computation of GSOP took 15 days. We cannot test the precise maximum size, when the computing times are such long. That is why we added another criterion, the result must be available in the time less then one hour. The sizes for GSOP in the first two columns are limited by this condition.

\subsection{Computing Time}
We tested also the computing times. There is one problem, the maximum sizes of Daoui and particularly Zhu algorithms are so low, that sufficient analysis of computing times is not possible. Finally, we tested these algorithms even if the error of orthogonality was higher than our threshold.

\begin{figure}[ht]
\centering
{\includegraphics[width=1\columnwidth,clip]{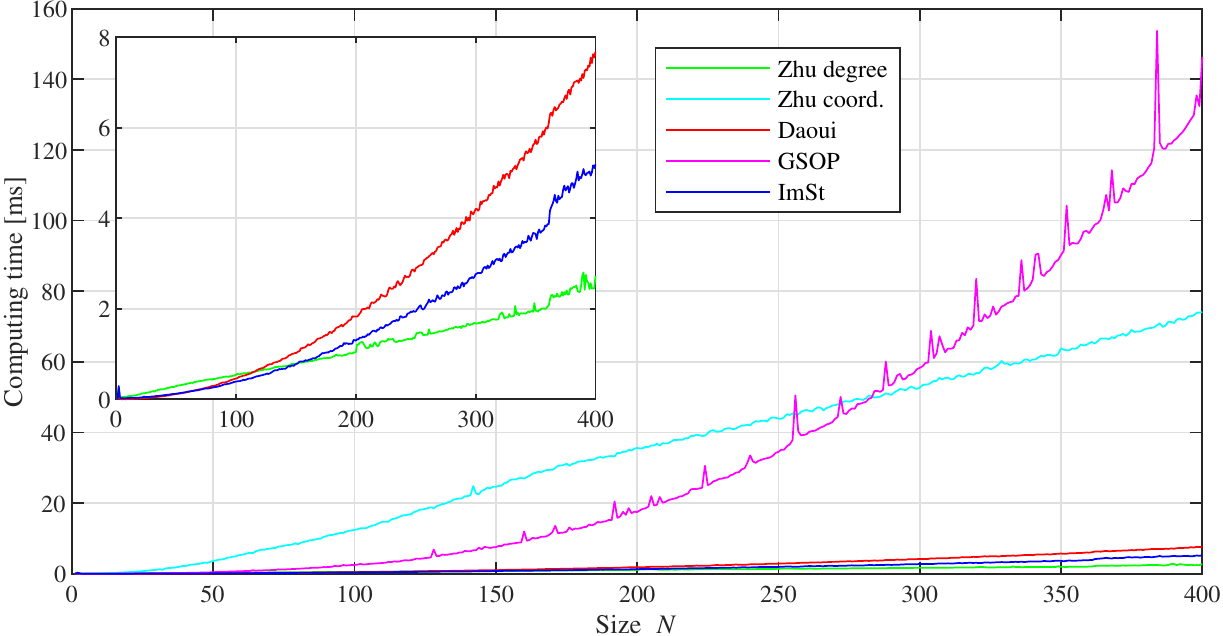}}
\caption{The computing times of the Racah polynomials for $n=N-1$, $a=\max(N/4,1)$, $b=a+N$, $\alpha=N/8$, and $\beta=N/16$.}
\label{fig:times}
\end{figure}

We choose these values of the parameters: $n=N-1$, $a=\max(N/4,1)$, $b=a+N$, $\alpha=N/8$, and $\beta=N/16$. We repeated each computation ten times, and took average time. The results are in \figurename{ \ref{fig:times}}.

The fastest algorithm is Zhu's recurrence over the degree, our algorithm ImSt is based on the similar principle, it is only a little bit slower. Daoui's algorithm is a little bit slower than ours and Zhu's recurrence over the coordinate is significantly slower, but it has still computing complexity $\mathcal{O}(N^2)$, only with worse constant. The computing complexity  $\mathcal{O}(N^3)$ of GSOP is clearly visible in the graph; from beginning, it is fast, but it cannot be used for high $N$.

\subsection{Restriction Error Analysis}
The distribution of moments is diverse from each other based on the discrete transforms \cite{CHP_2020}. To correctly reconstruct the signal information, the sequence of moments is important and should be recognized.
Therefore, the moment energy distribution of DRP is examined first; then the signal reconstruction analysis is performed. To acquire the distribution of moments, the procedure presented by Jian \cite{jain1989} is followed. The procedure is given in Algorithm \ref{alg3}.

\begin{algorithm}[ht]
	\caption{Find the moment order of DRP.}
	\label{alg3}
	{Input:}
	$\rho$=covariance coefficient \\
	{Output:}
	Order of DRP.
	\begin{algorithmic}[1]
		\State Generate the covariance matrix $ \Sigma $ with zero mean and length $ N $:
			\begin{equation}\label{Eq_cov_mat}
				\Sigma=\begin{bmatrix}
				1 		& 	\rho 	&  	\cdots 	& 	\rho^{N-1} \\
				\rho 	& 	1 		& 	     	&   \vdots 	\\
				\vdots 	& 	     	& \ddots 	& 	\rho \\
				\rho^{N-1} & 	\cdots 	& \rho 		& 	1 	
				\end{bmatrix}
			\end{equation}
	
		\State Transform the covariance matrix $ \Sigma $ into the domain of the discrete Racah moments ($M$) using
		\begin{equation}\label{Eq_CMM}
			M=R\times\Sigma\times R^T
		\end{equation}
	
		\State Find the diagonal coefficients $\sigma^2_{\ell}=M_{\ell\ell}$ of the discrete Racah moments $M$.
		
		\State Find the order of moments according to the values of the diagonal coefficients.
		\State \Return result
	\end{algorithmic}
\end{algorithm}

The covariance matrix $ \Sigma $ is used instead of an image. Then the matrix multiplication $ R\times\Sigma\times R^T $ can be used for moment computation, $ R $ is the matrix of Racah polynomials, $ R_{n,s-a}=\drpo{n}{s}$.
For the covariance coefficients, three values are used,
$ \rho=0.90 $, $ \rho=0.95 $, and $ \rho=0.98 $ with length $N=16$; then, the results are reported in \tablename{ \ref{tab:order}}.
From \tablename{ \ref{tab:order}}, it can observed that the maximum value of DRP is found at $ \ell=0 $ and the values are descendingly ordered. This declares that the DRP moment order used for signal reconstruction is $ n=0,1,\dots,N-1$.

\begin{table}[!httbp]
\caption{Transform coefficient values for different values of covariance coefficients.} 
\label{tab:order}
\resizebox{\textwidth}{!}{
\begin{tabular}{|c|c|c|c|c|c|c|c|c|c|c|c|c|}
\hline
 \multirow{5}{*}{$\ell$} & \multicolumn{4}{c|}{$\ $} & \multicolumn{4}{c|}{$\ $} & \multicolumn{4}{c|}{$\ $} \\[-1ex]
  & \multicolumn{4}{c|}{$\rho=0.9$} & \multicolumn{4}{c|}{$\rho=0.95$} & \multicolumn{4}{c|}{$\rho=0.98$} \\
 $\ $ & \multicolumn{4}{c|}{$\ $} & \multicolumn{4}{c|}{$\ $} & \multicolumn{4}{c|}{$\ $} \\[-2ex]
\cline{2-13}
  &&&&&&&&&&&&\\[-4.2ex]
& $a=0$ & $a=10$ & $a=30$ & $a=50$ & $a=0$ & $a=10$ & $a=30$ & $a=50$ & $a=0$ & $a=10$ & $a=30$ & $a=50$ \\
& $\alpha=a$ & $\alpha=a$ & $\alpha=a$ & $\alpha=a$ & $\alpha=a$ & $\alpha=a$ & $\alpha=a$ & $\alpha=a$ & $\alpha=a$ & $\alpha=a$ & $\alpha=a$ & $\alpha=a$ \\
& $\beta=0$ & $\beta=0$ & $\beta=0$ & $\beta=0$ & $\beta=0$ & $\beta=0$ & $\beta=0$ & $\beta=0$ & $\beta=0$ & $\beta=0$ & $\beta=0$ & $\beta=0$ \\
\hline
0 & 9.159 & 9.832 & 2.742 & 2.249 & 11.325 & 12.401 & 2.923 & 2.362 & 12.975 & 14.407 & 3.039 & 2.434 \\ \hline
1 & 2.912 & 2.856 & 2.388 & 2.045 & 2.232 & 1.907 & 2.500 & 2.125 & 1.527 & 0.916 & 2.568 & 2.174 \\ \hline
2 & 1.278 & 1.136 & 2.074 & 1.854 & 0.843 & 0.612 & 2.137 & 1.908 & 0.532 & 0.249 & 2.173 & 1.941 \\ \hline
3 & 0.702 & 0.591 & 1.794 & 1.675 & 0.440 & 0.300 & 1.822 & 1.708 & 0.272 & 0.120 & 1.836 & 1.727 \\ \hline
4 & 0.446 & 0.366 & 1.543 & 1.506 & 0.273 & 0.182 & 1.545 & 1.522 & 0.168 & 0.072 & 1.543 & 1.530 \\ \hline
5 & 0.311 & 0.252 & 1.313 & 1.345 & 0.188 & 0.124 & 1.298 & 1.347 & 0.115 & 0.049 & 1.286 & 1.347 \\ \hline
6 & 0.233 & 0.188 & 1.101 & 1.190 & 0.139 & 0.092 & 1.072 & 1.180 & 0.084 & 0.036 & 1.053 & 1.173 \\ \hline
7 & 0.183 & 0.147 & 0.900 & 1.037 & 0.109 & 0.072 & 0.862 & 1.018 & 0.065 & 0.028 & 0.838 & 1.005 \\ \hline
8 & 0.149 & 0.120 & 0.707 & 0.884 & 0.088 & 0.059 & 0.663 & 0.856 & 0.053 & 0.023 & 0.637 & 0.839 \\ \hline
9 & 0.125 & 0.101 & 0.523 & 0.725 & 0.074 & 0.050 & 0.477 & 0.691 & 0.044 & 0.020 & 0.450 & 0.671 \\ \hline
10 & 0.108 & 0.088 & 0.357 & 0.560 & 0.063 & 0.043 & 0.313 & 0.523 & 0.037 & 0.017 & 0.287 & 0.500 \\ \hline
11 & 0.095 & 0.077 & 0.224 & 0.396 & 0.055 & 0.038 & 0.183 & 0.357 & 0.032 & 0.015 & 0.159 & 0.334 \\ \hline
12 & 0.085 & 0.070 & 0.133 & 0.250 & 0.049 & 0.034 & 0.096 & 0.213 & 0.028 & 0.013 & 0.075 & 0.190 \\ \hline
13 & 0.077 & 0.063 & 0.083 & 0.142 & 0.044 & 0.031 & 0.051 & 0.108 & 0.025 & 0.012 & 0.031 & 0.088 \\ \hline
14 & 0.071 & 0.058 & 0.062 & 0.082 & 0.040 & 0.028 & 0.032 & 0.051 & 0.023 & 0.011 & 0.015 & 0.033 \\ \hline
15 & 0.066 & 0.054 & 0.055 & 0.058 & 0.037 & 0.027 & 0.027 & 0.030 & 0.021 & 0.010 & 0.011 & 0.014 \\ \hline
\end{tabular}
}
\end{table}

\begin{figure}[!htt]
\centering
\includegraphics[width = 0.75\textwidth]{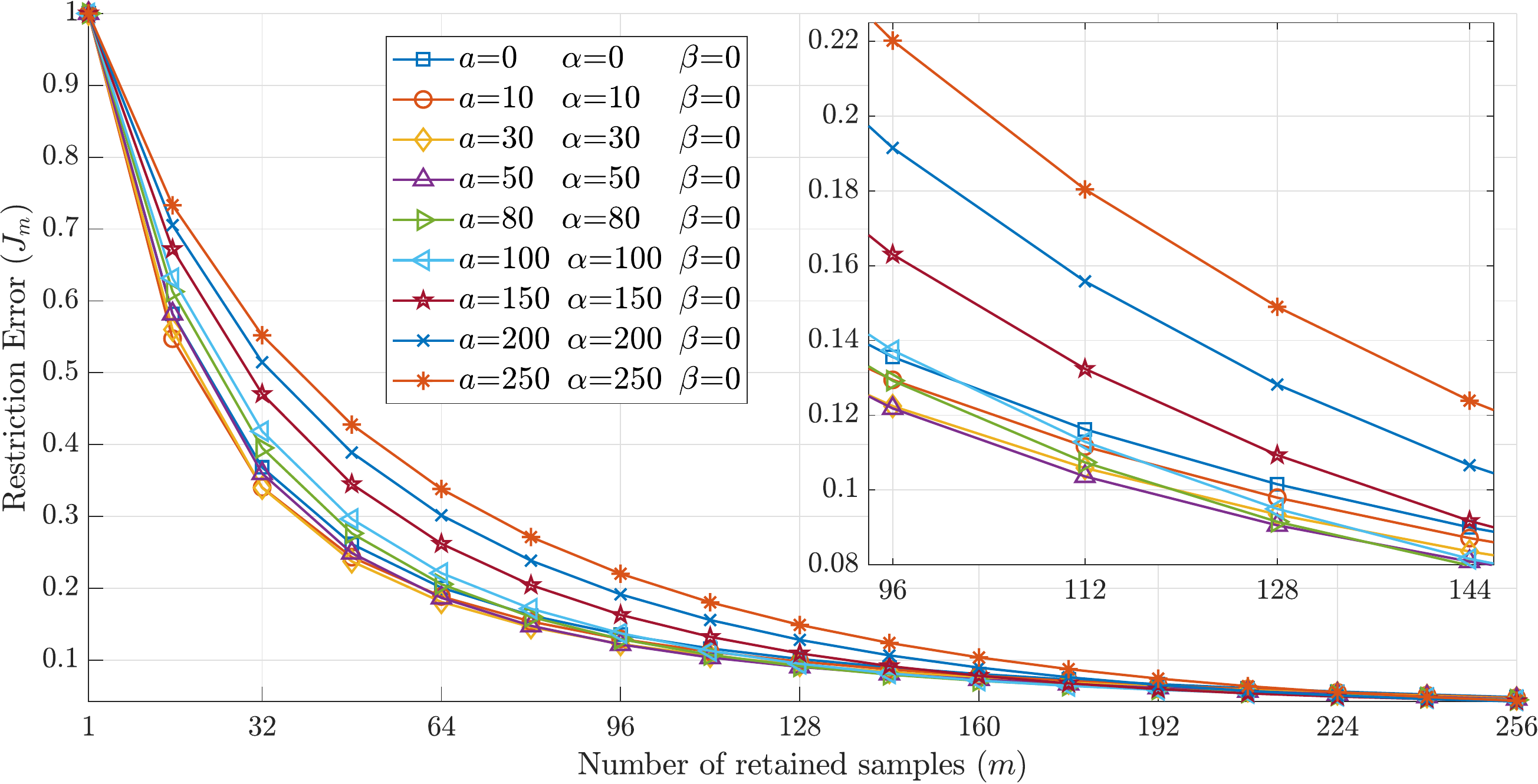}
\caption{Result of the restriction error for different values of Racah parameters values ($a=\alpha$ with $\beta=0$).}
\label{Fig:Rest1}
\end{figure}

The energy compaction property of the discrete transformation based on orthogonal polynomials is considered one of the important properties.
It is the fraction of the number of coefficients that reflect most of the signal energy to the total number of coefficients.
This characteristic is used to assess a DRP's ability to reconstruct a significant portion of the signal information from a very small number of moment coefficients.
To examine the impact of the DRP parameters $ a $, $ \alpha $ and $ \beta $ on the energy compaction, the restriction error, $ \mathcal{J} $, is used as follows \cite{jain1989}
\begin{equation}
	\label{Eq_Jm}
	\mathcal{J}_m=\frac{\sum\limits_{k=m}^{N-1}\sigma_k^2}{\sum\limits_{k=0}^{N-1}\sigma_k^2};\ \
	m=0,1,2,\dots,N-1 ,
\end{equation}
where $ \sigma_k^2 $ represents diagonal values of the transform coefficients ordered descendingly. In our case, the coefficients are already ordered, i.e. $k=\ell$.
\figurename{~\ref{Fig:Rest1}} shows the restriction error using the covariance coefficient $ \rho=0.95 $ with DRP parameters of $a=\alpha$ and $\beta=0$.
From \figurename{~\ref{Fig:Rest1}}, the DRP parameters affect the restriction error, which reveals that DRPs with parameters $ a=\alpha=30 $ and $ \beta=0 $ shows better energy compaction than other parameter values in the range of $ m<96 $. However, when $ a=\alpha=50 $ and $ \beta=0 $ presents better energy compaction compared to other DRP parameters in the range $ m>96 $.

\begin{figure}[!htt]
\centering
\includegraphics[width = 0.75\textwidth]{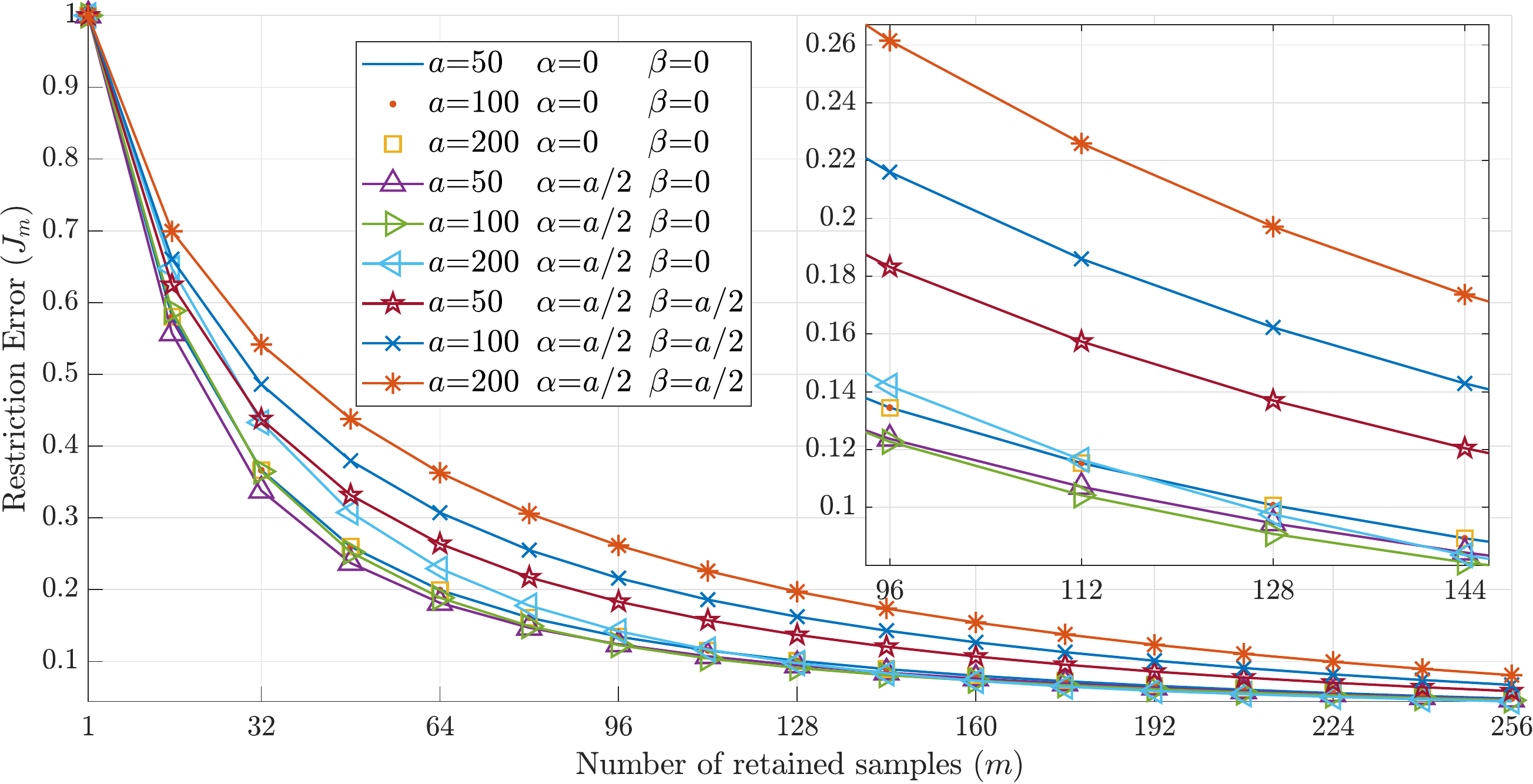}
\caption{Result of the restriction error for different values of Racah parameters values ($a$ and $\alpha=\{0,a/2\}$ with $\beta=\{0,a/2\}$).}
\label{Fig:Rest2}
\end{figure}

\figurename{~\ref{Fig:Rest2}} shows the restriction error of DRP with parameters of $a$, $\alpha=\{0,a/2\}$ and $\beta=\{0,a/2\}$. It can be observed from \figurename{~\ref{Fig:Rest2}} that DRPs with parameters $ a=50$, $\alpha=25 $ and $ \beta=0 $ shows better energy compaction than other parameter values in the range of $ m<96 $. However, when $ a=100$, $\alpha=50 $ and $ \beta=0 $ presents better energy compaction compared to other DRP parameters in the range $ m>96 $.

\begin{figure}[!htt]
\centering
\includegraphics[width = 0.75\textwidth]{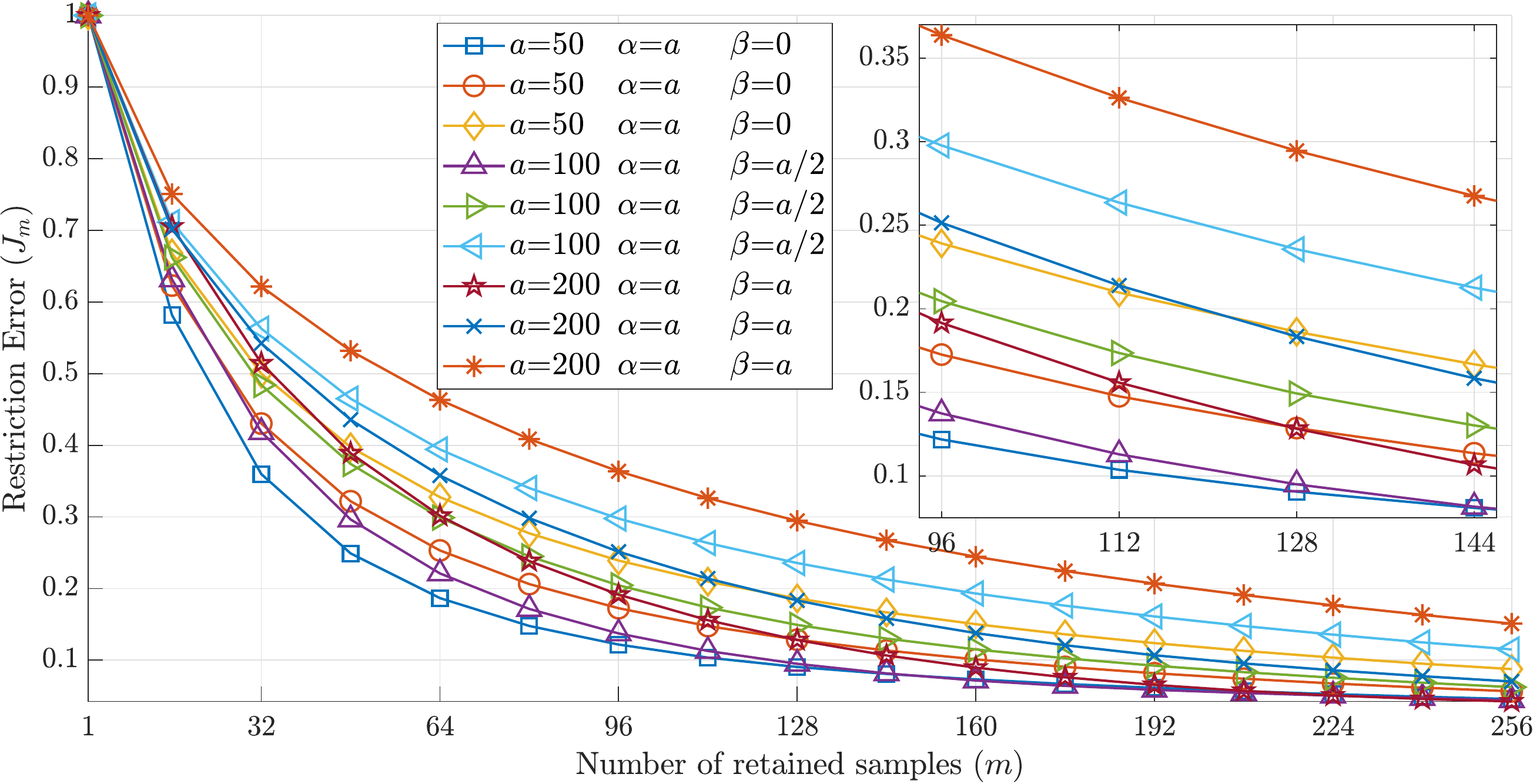}
\caption{Result of the restriction error for different values of Racah parameters values ($a$ and $\alpha=a$ with $\beta=\{0,a/2,a\}$).}
\label{Fig:Rest3}
\end{figure}

On the other hand, \figurename{~\ref{Fig:Rest3}} shows the restriction error for DRP with parameters of $a$, $\alpha=a$, and $\beta=\{0,a/2,a\} $. The best energy compaction in with this DRP parameters is $a=50$, $\alpha=50$, and $\beta=0 $ for the entire range of retained samples $ m $.

\subsection{Analysis of Reconstruction Error}

\begin{figure}[!ht]
\centering
\includegraphics[width = 0.50\textwidth]{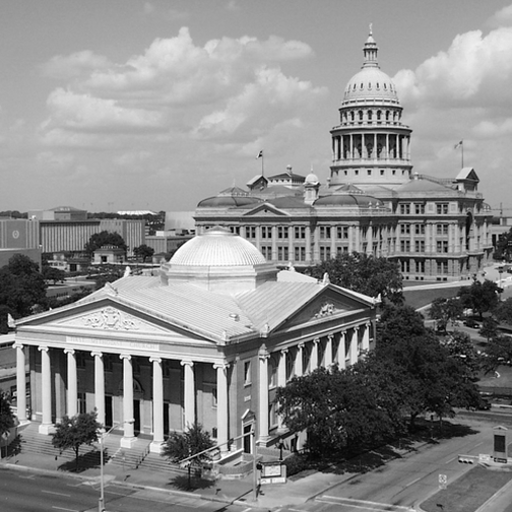}
\caption{Test image ``Church and Capitol''.}
\label{Fig:testimage}
\end{figure}

In this section, the reconstruction error analysis is carried out using real images.
The test image (``Church and Capitol''), shown in \figurename{ \ref{Fig:testimage}}, is taken from LIVE dataset \cite{LIVE}, \cite{LIVE_2006}, \cite{LIVE_2004}. The size of the test image is $512 \times 512$. Various values of DRP parameters ($a$, $\alpha$, $\beta$) were considered in the analysis. These parameters were considered in groups:
\begin{enumerate}
\item $a$ and $\alpha=\{0,a/2\}$ with $\beta=0$,
\item $a$ and $\alpha=\{0,a/2\}$ with $\beta=\{0,a/2\}$,
\item $a$ and $\alpha=a$ with $\beta=\{0,a/2,a\}$.
\end{enumerate}

DRPs ($R$) are generated first using the proposed sets of parameters. Then, DRMs ($M$) of the test image are computed.
After that, the image is reconstructed using the calculated moments and a finite number of moments. The normalized mean square error (NMSE), which compares the input image to the reconstructed version of the image, is then calculated. Thus, the NMSE is expressed as:
\begin{equation}
\label{NMSE}
    NMSE(I,I_r)=\frac{\sum\limits_{x,y} \left[I(x,y)-I_r(x,y)\right]^2}{\sum\limits_{x,y} I(x,y)^2},
\end{equation}
where $I$ and $I_r$ represent the original image and the reconstructed image, respectively.
NMSE is used as the reconstruction error.
From that the peek signal to noise ratio (PSNR)
\begin{equation}
\label{PSNR}
\begin{array}{l}
PSNR(I,I_r)=\\
\displaystyle
=10 \left(\log\left(\max\left(I(x,y)^2\right)\right)\!-\!\log\left(\frac{1}{N_1N_2}\sum\limits_{x,y} \left[I(x,y)\!-\!I_r(x,y)\right]^2\right)\right),
\end{array}
\end{equation}
where $N_1\times N_2$ is the size of the image $512\times 512$.

\begin{figure}[!ht]
\centering
\includegraphics[width = 0.75\textwidth]{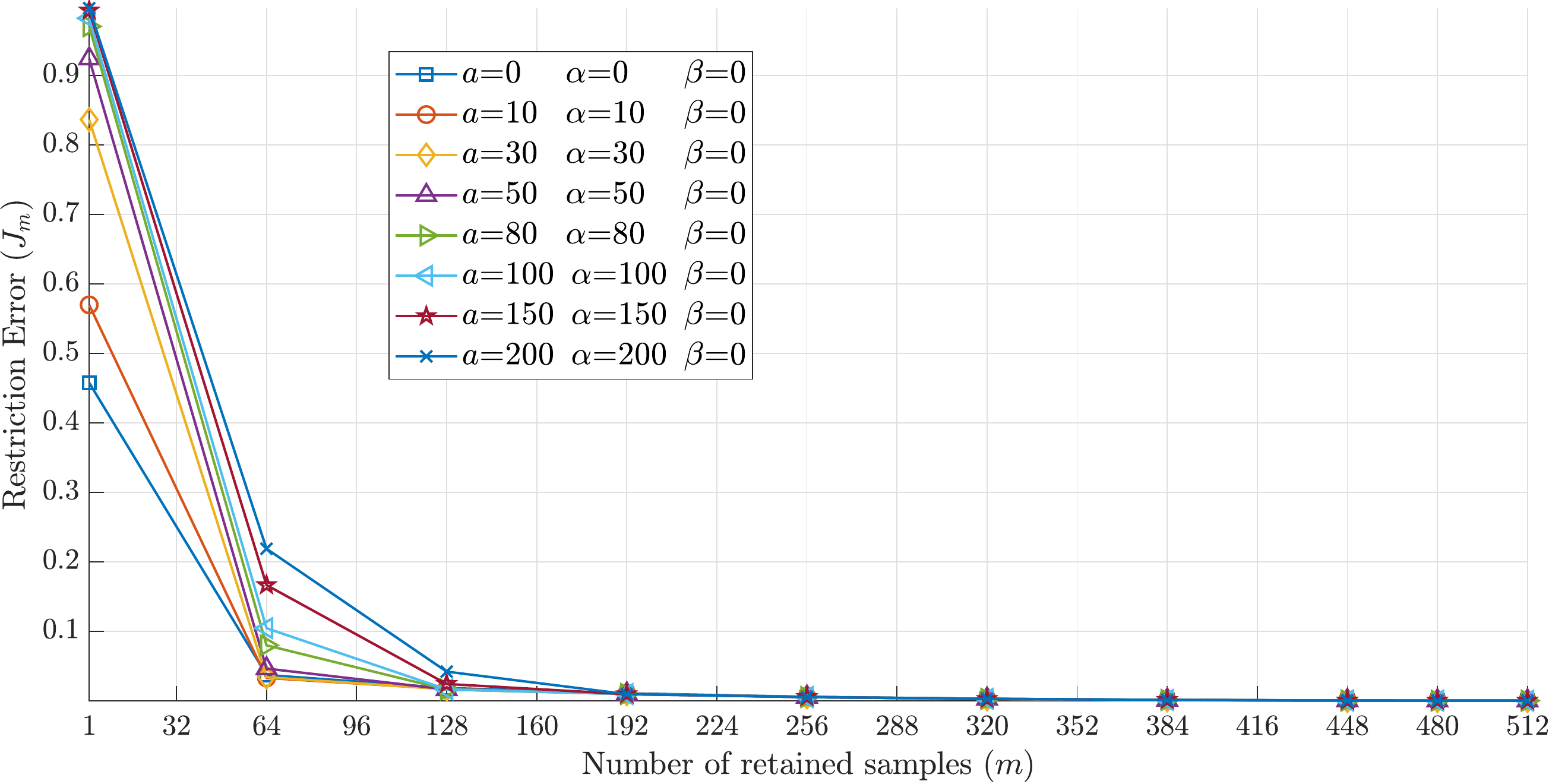}
\caption{Result of the reconstruction error using real image for different values of Racah parameters values ($a$ and $\alpha=a$ with $\beta=0$).}
\label{Fig:RecG1a}
\end{figure}

First, the reconstruction error analysis is carried out for $\alpha=a$ and $\beta=0$. The order of moments used to reconstruct the image is varied in the set $1, 32, 64, \dots,512$. The obtained results are depicted in \figurename{ \ref{Fig:RecG1a}}. The obtained results show that at moment order of 64, the best NMSE is occurred at DRP parameters of $a=10, \alpha=a, \beta=0$ with NMSE of 0.0328. The next three best NMSE are 0.0339, 0.0371, and 0.0463 for DRP parameters $a=30, \alpha=a, \beta=0$, $a=0, \alpha=a, \beta=0$, and $a=50, \alpha=a, \beta=0$, respectively. However, for moment order of 128, the NMSE are 0.0161, 0.0163, 0.01635, and 0.0165 for DRP parameters $a=80, \alpha=a, \beta=0$, $a=50, \alpha=a, \beta=0$, $a=100, \alpha=a, \beta=0$, and $a=30, \alpha=a, \beta=0$, respectively.

\begin{figure}[!ht]
\centering
\includegraphics[width = 0.99\textwidth]{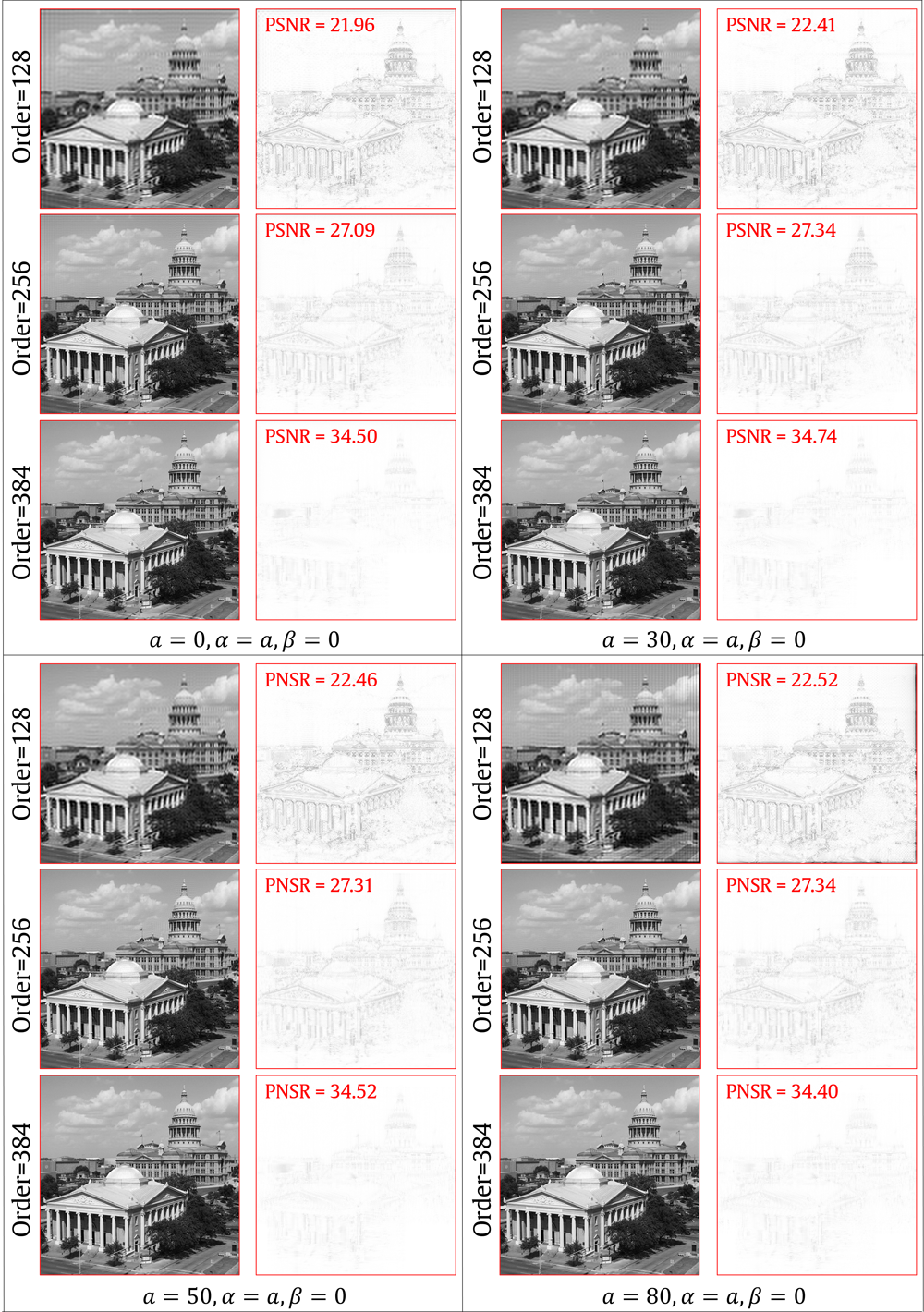}
\caption{Visual result of the reconstruction error using real image for different values of Racah parameters values ($a$ and $\alpha=(a)$ with $\beta=\{0,a/2\}$).}
\label{Fig:RecG1b}
\end{figure}

Moreover, for moment order of 256, the best NMSE is occurred at DRP parameters $a=10, \alpha=a, \beta=0$ with NMSE of 0.0052. For better inspection, Reconstruction error between the original and the reconstructed image is acquired and the PSNR is reported for different values of DRP parameters as shown \figurename{ \ref{Fig:RecG1b}}.

\begin{figure}[!ht]
\centering
\includegraphics[width = 0.75\textwidth]{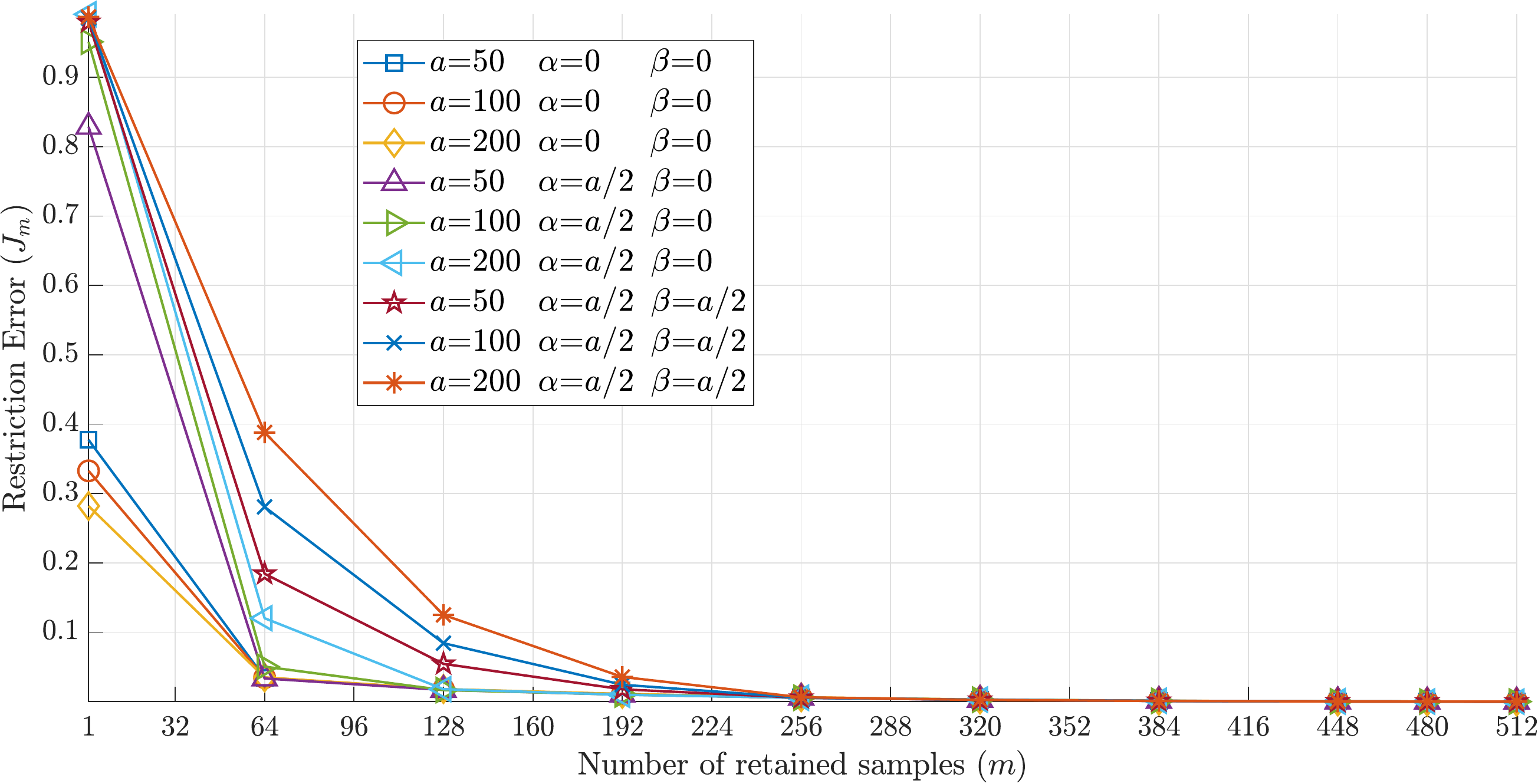}
\caption{Result of the reconstruction error using real image for different values of Racah parameters values ($a$ and $\alpha=\{0,a/2\}$ with $\beta=\{0,a/2\}$).}
\label{Fig:RecG2a}
\end{figure}

\begin{figure}[!ht]
\centering
\includegraphics[width = 0.99\textwidth]{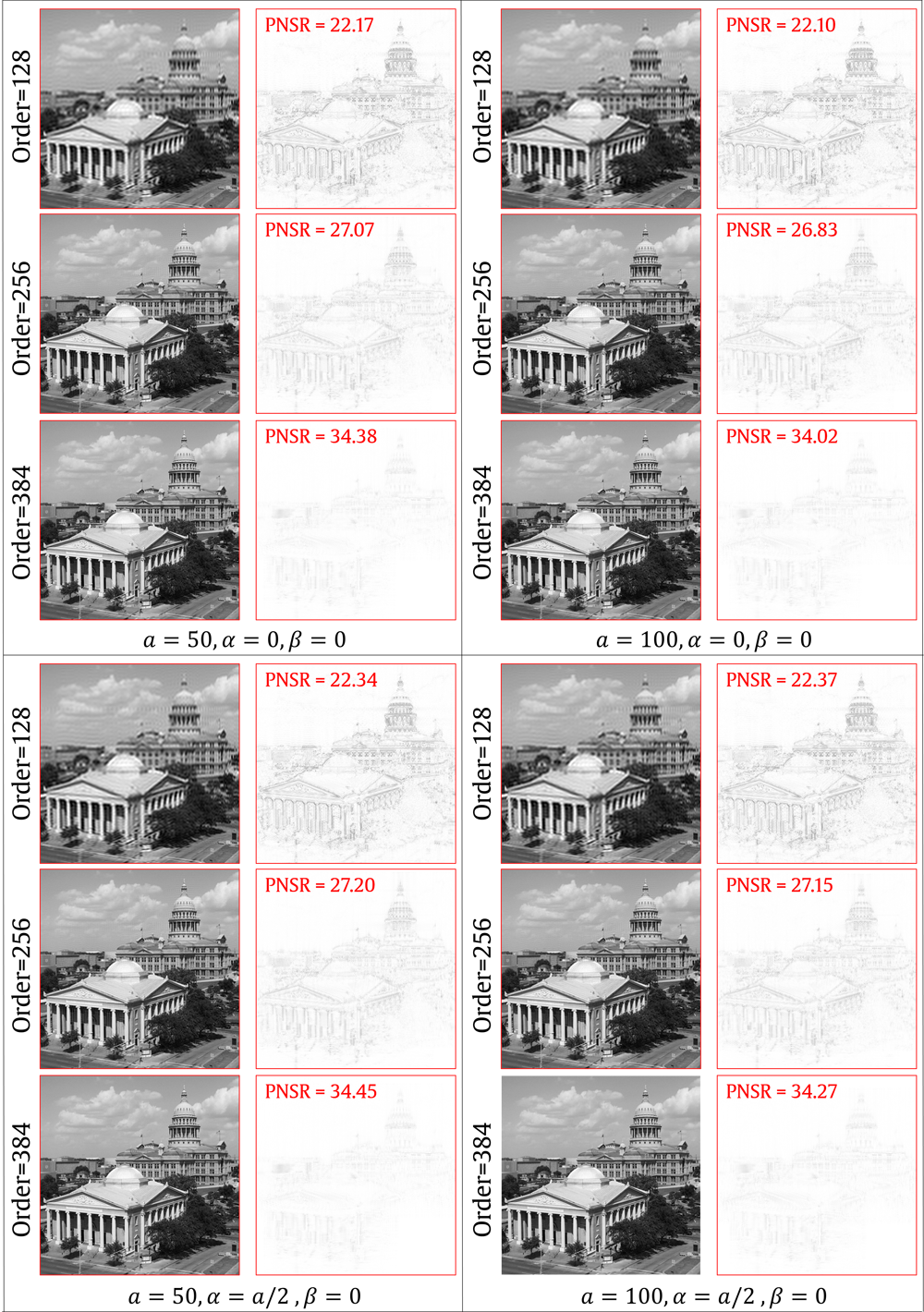}
\caption{Visual result of the reconstruction error using real image for different values of Racah parameters values ($a$ and $\alpha=\{0,a/2\}$ with $\beta=\{0,a/2\}$).}
\label{Fig:RecG2b}
\end{figure}

Second, the DRP parameter values in the range ($a$ and $\alpha=\{0,a/2\}$ with $\beta=\{0,a/2\}$) is used to perform the reconstruction error analysis.
The same moment orders of the first experiment is used in this experiment. \figurename{ \ref{Fig:RecG2a}} shows the obtained NMSE results of the second experiment. From \figurename{ \ref{Fig:RecG2a}}, the results demonstrate that at moment order of 64, the best NMSE is 0.0333 for DRP parameters of $a=50, \alpha=25, \beta=0$. The second best NMSE appears at $a=50, \alpha=0, \beta=0$ with NMSE of 0.0344; while the third best NMSE occurs at $a=100, \alpha=0, \beta=0$ with NMSE of 0.0347.
For moment order of 128, the best NMSE is 0.0166 for DRP parameters $a=100, \alpha=50, \beta=0$. In addition, the best NMSE,for moment order of 256, is occurred at DRP parameters $a=50, \alpha=25, \beta=0$ with NMSE of 0.00547.
For the sake of clarity, the visual reconstruction error between the original and the reconstructed image is acquired and the PSNR is reported for different values of DRP parameters as shown \figurename{ \ref{Fig:RecG2b}}.

\begin{figure}[!ht]
\centering
\includegraphics[width = 0.75\textwidth]{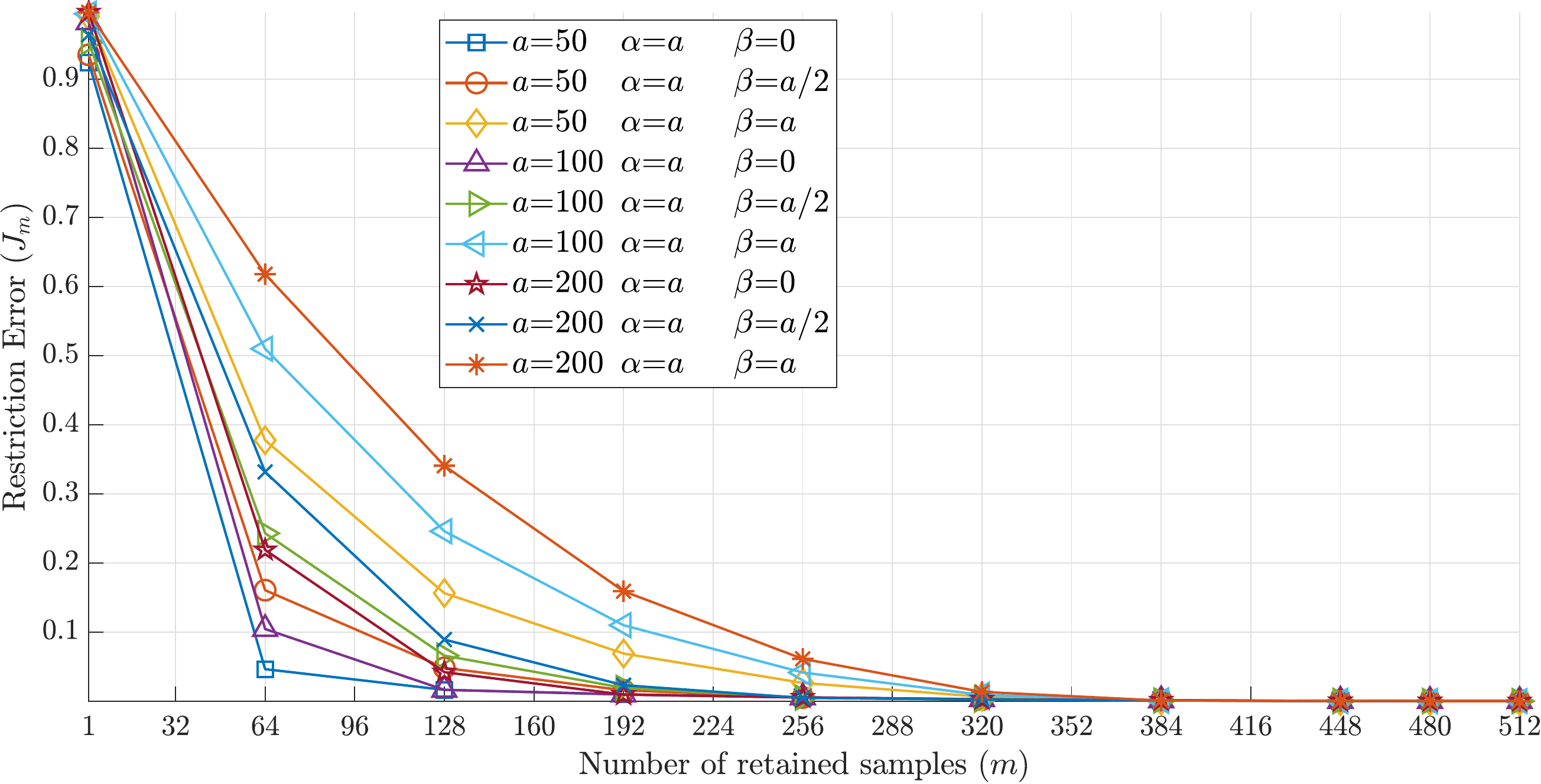}
\caption{Result of the reconstruction error using real image for different values of Racah parameters values ($a$ and $\alpha=a$ with $\beta=\{0,a/2,a\}$).}
\label{Fig:RecG3a}
\end{figure}

\begin{figure}[!ht]
\centering
\includegraphics[width = 0.99\textwidth]{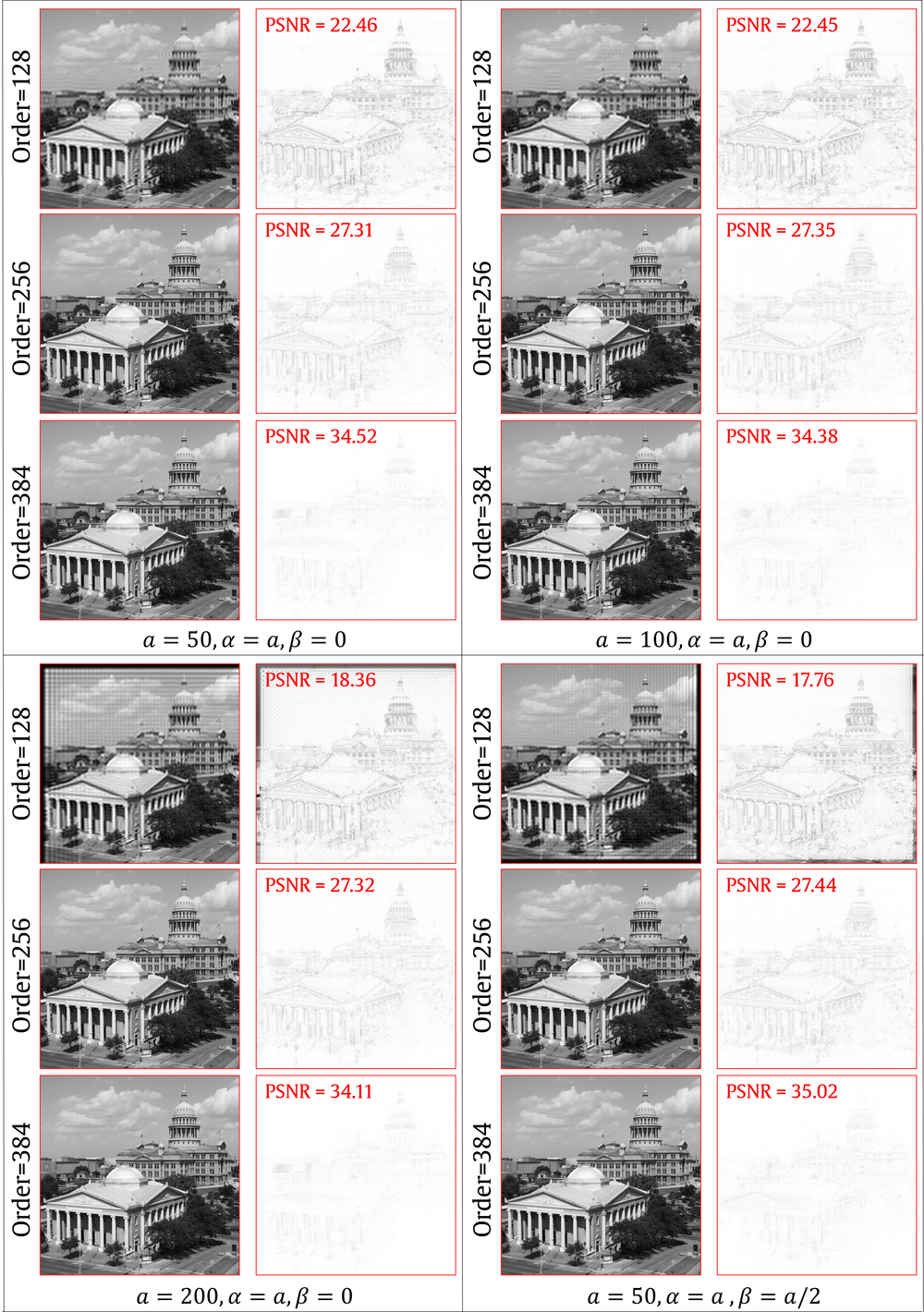}
\caption{Visual result of the reconstruction error using real image for different values of Racah parameters values ($a$ and $\alpha=a$ with $\beta=\{0,a/2,a\}$).}
\label{Fig:RecG3b}
\end{figure}

Finally, the DRP parameter values in the range ($a$ and $\alpha=a$ with $\beta=\{0,a/2,a\}$) is used to carry out the reconstruction error analysis.
\figurename{ \ref{Fig:RecG3a}} shows the reported NMSE results for this experiment. From \figurename{ \ref{Fig:RecG3a}}, the results demonstrate that at moment order of 64, the best NMSE is 0.0463 for DRP parameters of $a=50, \alpha=50, \beta=0$. The second best NMSE appears at $a=100, \alpha=100, \beta=0$ with NMSE of 0.104; while the third best NMSE occurs at $a=50, \alpha=50, \beta=25$ with NMSE of 0.161.
For moment order of 128, the best NMSE is 0.0163 for DRP parameters $a=50, \alpha=50, \beta=0$ and $a=100, \alpha=100, \beta=0$. In addition, the best NMSE, for moment order of 256, is occurred at DRP parameters $a=200, \alpha=200, \beta=100$ with NMSE of 0.00465.
The visual reconstruction error between the original and the reconstructed image is acquired and the PSNR is reported for different values of DRP parameters as shown \figurename{ \ref{Fig:RecG3b}}.

\section{Conclusion}
\label{Conc}
This paper proposed a new algorithm for computing the coefficient values of DRP. 
We use the logarithmic gamma function for computation the initial values. The utilization of the logarithmic gamma function empower the computation of the initial value for a wide range of DRP parameter values as well as large size of the polynomials. In addition, a new formula is used to compute the values of the initial sets based on the initial value.

The rest of DRP coefficients are computed by partitioning the DRP plane into four parts.
To compute the values in the four parts, the recurrence relation in the $x$ and $n$ directions are conjoined together. To clear out the propagation error, a stabilizing condition is forced. The performance of the proposed algorithm is tested against different values of DRP parameters. In addition, the proposed algorithm is compared with existing algorithms. These experiments show that the proposed algorithm reduced the computation cost compared to the existing algorithms. Moreover, the proposed algorithm is able to generate DRP for large sizes without propagation error. Finally, restriction error and reconstruction error analyses are performed to show the influence of the used parameter values.







\section*{Acknowledgments}

This work has been supported by the Czech Science Foundation (Grant No. GA21-03921S) and by the {\it Praemium Academiae}. We would also like to acknowledge University of Baghdad. 
for general and financial support.

\section*{Declarations}

We declare we have no conflict of interest.

\bibliographystyle{IEEEtran}


\end{document}